\documentclass[oneside, reqno, 11pt, a4paper]{amsart}

\usepackage{dutchcal}
\usepackage{microtype,times}
\AtBeginDocument{
  \DeclareSymbolFont{AMSb}{U}{msb}{m}{n}
  \DeclareSymbolFontAlphabet{\mathbb}{AMSb}
}
\DeclareFontFamily{U}{mathx}{\hyphenchar\font45}
\DeclareFontShape{U}{mathx}{m}{n}{<-> mathx10}{}
\DeclareSymbolFont{mathx}{U}{mathx}{m}{n}
\DeclareMathAccent{\widebar}{0}{mathx}{"73}

\usepackage[dvipsnames]{xcolor}
\usepackage{graphicx}
\usepackage{subfig}
\usepackage{svg}
\usepackage{tikz}
\tikzstyle{arrow} = [thick,->,>=stealth]
\usetikzlibrary{shapes.misc,matrix,fit,positioning,arrows.meta,decorations.pathreplacing,calc,shapes.geometric,arrows}
\tikzset{Matrix/.style={matrix of nodes, font=\footnotesize,text height=1pt, text depth=0.5pt, text width=8.5pt, align=center, column sep=0pt, row sep=0pt, nodes in empty cells}}

\usepackage[a4paper, pdftex, left=2cm, top=2cm, right=2cm, bottom=2cm]{geometry}
\usepackage[final]{pdfpages}
\usepackage{changepage}

\usepackage[foot]{amsaddr}

\usepackage{url}
\usepackage{hyperref}
\hypersetup{
    breaklinks=true,
    bookmarksopen=true,
    pdftitle={A novel space-time Discrete Empirical Interpolation Method}, 
    pdfauthor={Santiago Badia, Nicholas Mueller}, 
    pdfsubject={Nonaffinely parameterized partial differential equations, Space-time model order reduction, Reduced basis method, Discrete empirical interpolation, Hyper-reduction}, 
    pdfkeywords={complexity reduction, interpolation, parametric dependence}, 
    colorlinks=true,
    linkcolor=black,
    citecolor=blue,
    filecolor=black,
    urlcolor=blue
}

\newcommand{\quotes}[1]{``#1"}

\usepackage{csquotes}
\usepackage[backend=biber, maxbibnames=5, style=numeric-comp, isbn=false, bibencoding=utf8, safeinputenc, url=false, doi=true, giveninits=true]{biblatex}

\PassOptionsToPackage{normalem}{ulem}
\usepackage{ulem}
\providecolor{added}{rgb}{0,0,1}
\providecolor{deleted}{rgb}{1,0,0}

\usepackage{float}
\usepackage{framed}
\usepackage{verbatim}
\usepackage{fancyvrb}
\usepackage{booktabs}
\usepackage{colortbl}
\usepackage{multirow}
\usepackage{algorithm}
\usepackage{algpseudocode}
\makeatletter
\algdef{S}[IF]{IfNoThen}[1]{\algorithmicif\ #1}
\makeatother
\usepackage[inline]{enumitem}
\usepackage{siunitx}
\usepackage{mathtools}
\mathtoolsset{showonlyrefs}

\newtheorem{remark}{Remark}
\newtheorem{lemma}{Lemma}
\newtheorem{corollary}{Corollary}
\newtheorem{proposition}{Proposition}
\newtheorem{theorem}{Theorem}
\usepackage[breakable]{tcolorbox}
\tcbuselibrary{theorems}
\tcbuselibrary{skins}
\tcbset{
	commonstyle/.style={
		theorem style=plain,
		enhanced jigsaw,
		fonttitle=\bfseries,
		fontupper=\itshape,
		halign=justify,
		separator sign=:,
		description delimiters none,
		description font=\bfseries, 
		terminator sign={.\hspace{0.25em}},
		arc=0mm,outer arc=0mm,
		boxrule=0pt,toprule=0pt,bottomrule=0pt,leftrule=0pt,rightrule=0pt,
		titlerule=0pt,toptitle=0pt,bottomtitle=0pt,top=0pt,
		colback=white,coltitle=black,
		boxsep=0pt, bottom=0pt, left=0pt, 
	}
}
\newtcbtheorem[]{myproblem}{Problem}%
{center, commonstyle, fonttitle=\bfseries}{pb}
\newtcbtheorem[]{mydefinition}{Definition}
{center, commonstyle, fonttitle=\bfseries}{pb}
\newtcbtheorem[]{myassumption}{Assumption}
{center, commonstyle, fonttitle=\bfseries}{pb}

\usepackage{amsmath, amsfonts, amssymb, amscd, bm, mathtools}

\renewcommand{\vec}[1]{\underline{#1}}

\newcommand{\N}{\mathbb{N}}
\newcommand{\R}{\mathbb{R}}
\newcommand{\abs}[1]{\left| #1 \right|}
\newcommand{\norm}[1]{\vert\vert #1 \vert\vert}

\newcommand{\spod}[1]{\mathrm{TPOD}_{s, \hat s}\left( #1 \right)}
\newcommand{\tpod}[1]{\mathrm{TPOD}_{t, \hat t}\left( #1 \right)}
\newcommand{\roundbrackets}[1]{\left( #1 \right)}
\newcommand{\squarebrackets}[1]{\left[ #1 \right]}

\DeclareMathOperator*{\argmax}{\arg\max}
\newcommand{\zero}[1]{\bm{0}_{#1}}

\usepackage{acronym}
\acrodef{pde}[PDE]{partial differential equation}
\acrodef{fe}[FE]{finite element}
\acrodef{fem}[FEM]{finite element method}
\acrodef{dof}[DOF]{degree of freedom}
\acrodefplural{dof}[DOFs]{degrees of freedom}
\acrodef{be}[BE]{Backward Euler}
\acrodef{hf}[HF]{high fidelity}
\acrodef{fom}[FOM]{full order model}
\acrodef{lhs}[LHS]{left hand side}
\acrodef{rhs}[RHS]{right hand side}
\acrodef{rom}[ROM]{reduced order model}
\acrodef{rb}[RB]{reduced basis}
\acrodefplural{rb}[RBs]{reduced bases}
\acrodef{svd}[SVD]{singular value decomposition}
\acrodef{sthosvd}[ST-HOSVD]{sequentially truncated high order singular value decomposition}
\acrodef{tpod}[TPOD]{truncated proper orthogonal decomposition}
\acrodef{strb}[ST-RB]{space-time reduced basis}
\acrodef{mdeim}[MDEIM]{empirical interpolation method in matrix form}
\acrodef{stmdeimrb}[ST-MDEIM-RB]{space-time MDEIM-RB}
\acrodef{stdmdeim}[STD-MDEIM]{standard MDEIM}
\acrodef{stmdeim}[ST-MDEIM]{space-time MDEIM}
\acrodef{funmdeim}[FUN-MDEIM]{functional version of MDEIM}
\acrodef{stfunmdeim}[STFUN-MDEIM]{space-time functional MDEIM}

\graphicspath{{img/}}

\newcommand{\dbhat}[1]{\widehat{\raisebox{0pt}[0.9\height]{$\widehat{#1}$}}}

\title[Model order reduction with novel discrete empirical interpolation methods in space-time]{Model order reduction with novel discrete empirical interpolation methods in space-time}
\date{\today}

\address{$^\dagger$School of mathematics\\Monash university\\Clayton\\Victoria 3800\\Australia}
\address{$^\ast$Centre Internacional de M\`etodes Num\`erics a l'Enginyeria, Campus Nord, 08034, Barcelona, Spain.}
\author[N. Mueller]{Nicholas Mueller$^\dagger$}
\email{nicholas.mueller@monash.edu}
\author[S. Badia]{Santiago Badia$^{\dagger\ast}$}
\email{santiago.badia@monash.edu}

\begin{document}

\begin{abstract}
  This work proposes novel techniques for the efficient numerical simulation of parameterized, unsteady partial differential equations. Projection-based reduced order models (ROMs) such as the reduced basis method employ a (Petrov-)Galerkin projection onto a linear low-dimensional subspace. In unsteady applications, space-time reduced basis (ST-RB) methods have been developed to achieve a dimension reduction both in space and time, eliminating the computational burden of time marching schemes. However, nonaffine parameterizations dilute any computational speedup achievable by traditional ROMs. Computational efficiency can be recovered by linearizing the nonaffine operators via hyper-reduction, such as the empirical interpolation method in matrix form. In this work, we implement new hyper-reduction techniques explicitly tailored to deal with unsteady problems and embed them in a ST-RB framework. For each of the proposed methods, we develop a posteriori error bounds. We run numerical tests to compare the performance of the proposed ROMs against high-fidelity simulations, in which we combine the finite element method for space discretization on 3D geometries and the Backward Euler time integrator. In particular, we consider a heat equation and an unsteady Stokes equation. The numerical experiments demonstrate the accuracy and computational efficiency our methods retain with respect to the high-fidelity simulations.
\end{abstract}

\maketitle

\section{Introduction}
\label{sec: introduction}
Conventional high-fidelity approximations of phenomena governed by \acp{pde} generally involve fine space and time discretizations, requiring large computational resources. When the \ac{pde} is parameterized, the parameter-to-solution map is a manifold embedded in a high-dimension space. Instead of performing the \ac{hf} simulation every time we want to get the solution for a given parameter, one can design a surrogate model for the parameter-to-solution map.  Projection-based \acp{rom}, such as the \ac{rb} method, create surrogate models by approximating the solution manifold by a lower-dimensional vector space. These surrogates provide accurate solutions in a much shorter time and with much fewer computational resources. 

Traditionally, these methods only aim to reduce the spatial complexity of the \ac{fom} by employing (Petrov-)Galerkin projections on the ordinary differential equation (ODE) resulting from the space-only discretization, while using the same \ac{hf} time marching scheme for the integration in time. Space-time \acp{rom} can be applied to the fully discrete system of \ac{fom} equations to retain a reduction along the temporal and spatial dimensions. In particular, the \ac{strb} method initially builds a spatial and a temporal basis spanning a low-dimensional subspace, using either a greedy algorithm or a \ac{tpod} approach. Then, under a suitable norm, it minimizes the \ac{fom} residual on such subspace. 

A comprehensive literature review of methods to reduce temporal complexity can be found in \cite{choi2019space}. The methods presented in \cite{urban2012new, urban2014improved, yano2014space, yano2014space2} are characterized by a reduction of both spatial and temporal dimensions. However, they exhibit some relevant drawbacks, the most salient ones being the need for (uncommon) \ac{fe} discretizations of the time domain and the absence of hyper-reduction techniques to handle non-linearities efficiently. The authors of \cite{choi2021space, kim2021efficient, Tenderini2022Mueller} successfully achieve a temporal complexity reduction by implementing \ac{strb} in the case of both 2D and 3D problems, demonstrating promising results both in terms of accuracy and computational efficiency. However, the applications discussed in these works are restricted to simple parameter-to-solution maps, i.e., affine dependence of operators on parameters.  

In general, the performance of linear \acp{rom} such as \ac{strb} depends on whether the \ac{fom} is reducible or not, which in turn is determined by two main factors:
\begin{itemize}
    \item There exists a low-dimensional spatio-temporal manifold that accurately approximates the manifold of the \ac{fom} solutions.
    \item The implementation of the method admits an efficient splitting into an \emph{offline phase}, (where all the expensive, parameter-independent operations are executed once and for all) and an \emph{online phase} (where the inexpensive parameter-dependent computations are performed for every parameter of interest). 
\end{itemize}
The first condition is generally met when dealing with elliptic or dissipative time-dependent applications whose parameter-to-solution map is sufficiently regular \cite{Unger_2019,quarteroni2015reduced}. The second condition stems from the parametric complexity of the problem. It is satisfied when the parameterizations of the \ac{fom} are affine, i.e. they can be written as a linear combination of terms. When this does not occur, we can use a hyper-reduction technique to approximate the nonaffinities as a sum of a small number of affine quantities. 

Hyper-reduction strategies are divided into two categories: collocation procedures \cite{Legresley2006Alonso,Astrid2008Weiland,RYCKELYNCK2005346}, such as the Gauss-Newton approximated tensors (GNAT) method \cite{CARLBERG2013623}, which rely on sampling procedures applied to the nonaffine function; and function-reconstruction approaches, such as the empirical interpolation method (EIM) \cite{Everson:95,BARRAULT2004667,https://doi.org/10.1002/fld.2712,refId0,https://doi.org/10.1002/nme.2309}, its discrete variant (DEIM) \cite{chaturantabut2010nonlinear,https://doi.org/10.1002/nme.3050,doi:10.1137/120899042,FOSASDEPANDO2016194}, and the \ac{mdeim} \cite{doi:10.1137/140959602,doi:10.1137/130932715,https://doi.org/10.1002/fld.4260}. The latter is a purely \emph{algebraic} interpolation approach that was proposed to approximate parameterized bilinear operators. Although \ac{mdeim} has been successfully applied to unsteady problems \cite{NEGRI2015431,ELZOHERY2022109396}, it was not originally conceived to reduce the temporal complexity of time-dependent bilinear operators. Additionally, the cost of its offline phase explodes when dealing with unsteady problems since an evaluation of the algebraic operator is required at every time step. 

We list hereunder the contributions of our work:
\begin{itemize}
    \item We extend the implementation of \ac{mdeim} to compress bilinear operators in space and time. In particular, we approximate them as a linear combination of \ac{strb} basis vectors by coefficients satisfying an interpolation condition on a set of space-time indices chosen via a greedy algorithm. The same idea can be used to approximate the \ac{rhs} of the problem.
    \item We propose a \emph{functional} alternative to the algebraic implementation of \ac{mdeim} to reduce the offline cost. Here, we first approximate the fields in the bilinear forms that depend nonlinearly on the temporal variable and the parameters (e.g. physical coefficients), similarly to what is done in \cite{https://doi.org/10.1002/fld.2712,doi:10.1137/10081157X,LASSILA20101583,Antil2014}. Then, the bilinear operators obtained by substituting the nonaffine fields with their compressed versions are approximated with the standard algebraic \ac{mdeim}. The computational cost of this step is greatly reduced thanks to the preliminary compression of the fields.
    \item We integrate the proposed hyper-reductions with a \ac{strb} state approximation, which results in a \ac{stmdeimrb} method. We derive \emph{a posteriori} error bounds for each of the \ac{stmdeimrb} approaches, which we show decay according to a user-defined tolerance that determines the accuracy of the methods, under suitable well-posedness and reducibility assumptions.
    \item We investigate the numerical properties of the proposed \ac{stmdeimrb} approaches to the heat equation with parameterized thermal diffusivity and boundary values and the unsteady Stokes equations with a paramaterized and time-dependent fluid viscosity and inflow rate. Both tests are conducted on 3D geometries.
\end{itemize}

This article is organized as follows. In Sect.~\ref{sec: notation}, we provide the notation we employ throughout our work. In Sect.~\ref{sec: section2}, we discuss the nature and the implementation of the \ac{strb} approaches. In Sect.~\ref{sec: section3}, we review the current literature on \ac{mdeim} for steady problems, and we propose our novel \ac{mdeim} strategies for unsteady applications. In Sect.~\ref{sec: section4}, we present the ST-MDEIM-RB  methods and prove \emph{a posteriori} error estimates. The derivations in Sections \ref{sec: section2}-\ref{sec: section4} are shown using the heat equation as model problem. We consider the Stokes equations in Sec.~\ref{sec: section5}. In Sect.~\ref{sec: section6}, we discuss the numerical results obtained when applying ST-MDEIM-RB to our test cases. Finally, in Sect.~\ref{sec: section7} we draw some conclusions and possible extensions.

\section{Hypermatrices and basic notation}
\label{sec: notation}
In this work, we manipulate 3-dimensional arrays, which we call hypermatrices. Given a hypermatrix, we label its axes with specific subscripts (kept separate by commas). In particular, we use the subindices $s$, $t$ and $\bm{\mu}$ to indicate the spatial, temporal and parametric axes, respectively. For example, $\bm{U}_{s,t,\bm{\mu}} \in \R^{N_s \times N_t \times N_{\bm{\mu}}}$ can be interpreted as an array of $N_{\bm{\mu}}$ \ac{hf} (snapshot) matrices of size $N_s \times N_t$, where $N_s$, $N_t$ and $N_{\bm{\mu}}$ represent the number of \ac{hf} spatial \acp{dof}, \ac{hf} time steps, and parameters. Our algorithms extensively employ the \ac{tpod} operation to compress hypermatrices in specific indices. We use the subscripts $\hat s$ and $\hat t$ to indicate a reduction in the spatial and temporal dimensions, respectively, and we denote the number of reduced spatial and temporal \acp{dof} as $n_s$ and $n_t$.  

We now define two basic operations on hypermatrices, i.e. \emph{permutation of axes} and \emph{flattening}. The former is a generalization of matrix transposition in multiple dimensions. For example, $\bm{U}_{t,s,\bm{\mu}} \in \R^{N_t \times N_s \times N_{\bm{\mu}}}$ is the result of permuting space and time axes of $\bm{U}_{s,t,\bm{\mu}}$, and is marked by a different ordering of the subindices. The latter is equivalent to the operation of \emph{merging axes}, e.g., $\bm{U}_{s,t\bm{\mu}} \in \R^{N_s \times N_t N_{\bm{\mu}}}$ is the result of flattening the time and parameter indices of $\bm{U}_{s,t,\bm{\mu}}$ into a single index, ending up with a matrix. Since these two operations are isomorphisms, we do not need to distinguish between the hypermatrix before and after the permutation/flattening. We use the standard notation $\cong$ to denote congruence by isometry of two objects. For example, $\bm{U}_{t,s \bm{\mu}} \cong \bm{U}_{s,t,\bm{\mu}}$, as they provide the same information while being indexed differently (due to flattening and permutation of indices). For the product of hypermatrices, we use the same notation as for matrices, e.g., $\bm{U}_{s,\hat s} \bm{V}_{\hat s} \in \R^{N_s}$.  

The computation of a hypermatrix \ac{tpod} is a generalization of a matrix standard \ac{tpod}. We define the space-axis \ac{tpod} as:  
\begin{equation}
    \label{eq: truncated POD space}
    \bm{\Phi}_{s, \hat s} =  \underset{\bm{Z}_{s,\hat s} \in \R^{N_s \times n_s} : \ \bm{Z}_{\hat s, s}\bm{Z}_{s,\hat s} = \bm{I}_{\hat s, \hat s}}{\mathrm{arg min}} \|\bm{I}_{s,t \bm{\mu}} -  \bm{Z}_{s, \hat s} \bm{Z}_{\hat s, s} \bm{U}_{s,t\bm{\mu}} \|_F,
\end{equation}
where $\bm{I}$ stands for an identity matrix (with axes determined by the subindices), $\norm{\cdot}_F$ indicates the Frobenius norm of a matrix, and $n_s$ is the rank of the truncation, defined as in \eqref{eq: truncated POD}. We are flattening the time and parameter axes and compressing the $N_t N_{\bm{\mu}}$ space snapshots into $n_s$ space basis vectors. The result is a matrix $\bm{\Phi}_{s, \hat s}$ that collects the reduced space basis vectors in the columns. We interpret such a matrix as a basis spanning an $n_s$-dimensional spatial subspace of $\R^{N_s}$ on which we approximate the spatial evolution of the \ac{hf} solution. The matrix $\bm{\Phi}_{\hat s, s}$ is the transpose of $\bm{\Phi}_{s, \hat s}$ and is used to spatially compress the snapshots by projecting them onto $\mathrm{range}(\bm{\Phi}_{\hat s, s})$. Using our notation, this corresponds to computing $\widehat{\bm{U}}_{\hat s,t \bm{\mu}} =  \bm{\Phi}_{\hat s,s} \bm{U}_{s,t \bm{\mu}}$, where $\widehat{\cdot}$ denotes a compressed quantity, i.e. it can be written as a linear combination of the \ac{rb} vectors. Note that $\widehat{\bm{U}}_{s,t \bm{\mu}} =  \bm{\Phi}_{s,\hat s} \widehat{\bm{U}}_{\hat s,t \bm{\mu}}$ is still in the reduced space, even though written in the \ac{hf} basis.

Different options are available to achieve a reduction in time. In this work, we choose the \ac{sthosvd} approach defined in \cite{carlberg2017galerkin}. We first compute the space-axis \ac{tpod}, and then we compute the time-axis \ac{tpod} of the reduced space snapshots as:
\begin{equation}
    \label{eq: truncated POD time}
    \bm{\Phi}_{t, \hat t} =  \underset{\bm{Z}_{t,\hat t} \in \R^{N_t \times n_t} : \ \bm{Z}_{\hat t,t}\bm{Z}_{t,\hat t} = \bm{I}_{\hat t, \hat t}}{\mathrm{arg min}} \|\bm{I}_{t,\hat s \bm{\mu}} -  \bm{Z}_{t, \hat t} \bm{Z}_{\hat t, t} \widehat{\bm{U}}_{t,\hat s \bm{\mu}} \|_F,
\end{equation}
Here, we are compressing the $n_s N_{\bm{\mu}}$ space-compressed and parameter snapshots into $n_t$ time basis vectors. The time-axis \ac{tpod} on the already compressed hypermatrix in space aims to reduce the computational cost of the operation. Similarly to the spatial basis, we interpret $\bm{\Phi}_{t, \hat t}$ as a basis spanning an $n_t$-dimensional temporal subspace of $\R^{N_t}$ on which we approximate the temporal evolution of the \ac{hf} solution.

We can now define the space-time \ac{rb} as 
\begin{equation}
    \label{eq: truncated POD space-time}
\bm{\Phi}_{st,\hat s \hat t} = \bm{\Phi}_{s,\hat s} \otimes \bm{\Phi}_{t,\hat t}, 
\end{equation}  
where $\otimes$ denotes the standard Kronecker product between matrices. We aim to approximate the evolution of the \ac{hf} solution in space-time on the subspace spanned by the $n_{st} = n_sn_t$ columns of $\bm{\Phi}_{st,\hat s \hat t}$. We denote the space-time compression of $\bm{U}_{s,t}$ as $\widehat{\bm{U}}_{\hat s, \hat t} \cong \bm{\Phi}_{\hat s \hat t, st} \bm{U}_{st}$, whereas $\widehat{\bm{U}}_{s,t} \cong \bm{\Phi}_{st, \hat s \hat t} \widehat{\bm{U}}_{\hat s \hat t}$ is the expansion of the compressed quantity on the \ac{hf} \acp{dof}. Lastly, $A \lesssim B$ (resp. $A \gtrsim B$, $A \eqsim B$) denotes $A \leq c B$ (resp. $A \geq c B$, $A = c B$ ) for a constant $c$; and we use $\N(M)$ to indicate $\{  1,\ldots,M\}$, with $M \in \N_+$.

\section{Space-time reduced basis methods for a parameterized heat equation}
\label{sec: section2}
We commence this section by introducing the \ac{fom} associated with a parameter-dependent heat equation. Consequently, we present in detail the offline and online phases of the \ac{strb} method applied to the \ac{fom}.

\subsection{Full order model} 
\label{subs:full order model}
We consider an open, polyhedral Lipschitz domain $\Omega \subset \R^d$ with boundary $\partial\Omega$, and a temporal domain $[0,T] \subset \R_+ \cup \{0\}$. We also consider a parameter space $\mathcal{D} \subset \R^p$. Given $\bm{\mu} \in \mathcal{D}$, a generic heat equation defined on the space-time domain $\Omega \times [0,T]$ reads:
\begin{equation} 
\label{eq: strong form heat equation}
\begin{cases}
\frac{\partial u^{\bm{\mu}}}{\partial t} - \bm{\nabla} \cdot (\alpha^{\bm{\mu}} \bm{\nabla} u^{\bm{\mu}}) = f^{\bm{\mu}}  & \text{in} \ \Omega \times (0,T], \\
u^{\bm{\mu}} = g^{\bm{\mu}}  & \text{on} \ \Gamma_D \times (0,T], \\
\alpha^{\bm{\mu}} \vec{n} \cdot \bm{\nabla} u^{\bm{\mu}}  = h^{\bm{\mu}} & \text{on} \ \Gamma_N \times (0,T], \\
u^{\bm{\mu}} = u_0^{\bm{\mu}} & \text{in} \ \Omega \times \{0\},
\end{cases} 
\end{equation}
where $u^{\bm{\mu}}: \Omega \times [0,T] \to \R$ is the unknown temperature (state variable), $\alpha^{\bm{\mu}}: \Omega \times [0,T] \to \R_{+}$ is the heat diffusivity coefficient, $f^{\bm{\mu}}: \Omega \times [0,T] \to \R$ is the forcing term, $g^{\bm{\mu}}: \Gamma_D \times [0,T] \to \R$ and $h^{\bm{\mu}}: \Gamma_N \times [0,T] \to \R$ are the Dirichlet and Neumann boundary data, respectively, $u_0^{\bm{\mu}}: \Omega \to \R$ is the initial condition, $\vec{n}$ is the outward unit normal vector to $\partial\Omega$, and $\Gamma_D$ and $\Gamma_N$ are respectively the Dirichlet and Neumann boundaries (such that $\{\Gamma_D, \Gamma_N\}$ forms a partition of $\partial\Omega$). The superindex $(\cdot)^{\bm{\mu}}$ indicates the dependence of a variable $(\cdot)$ upon $\bm{\mu}$.

Let us introduce the Hilbert spaces
\begin{equation}
    \mathcal{V} = H^1(\Omega); \qquad \mathcal{V}^0_{\Gamma_D}  
    = \left\{ v \in H^1(\Omega) \ \ : \   v=0 \  \text{on} \ \Gamma_D\right\}.
\end{equation}
Let us consider a conforming quasi-uniform partition $\mathcal{T}_h$ of $\Omega$ with meshsize $h$ and a uniform partition of the time domain $(0,T)$ with time-step size $\delta = T/N_t$; we define the sequence of time steps $\{t_n\}_{n=0}^{N_t}$ with $t_{n} = n \delta$. For the spatial approximation of \eqref{eq: strong form heat equation}, we consider a conforming \ac{fe} space $\mathcal{V}_h \subset \mathcal{V} $ on $\mathcal{T}_h$ and $\mathcal{V}_h^0 = \mathcal{V}_h \cap \mathcal{V}^0_{\Gamma_D}$. For the time discretization, we consider the \ac{be} scheme for simplicity in the exposition. To strongly impose non-homogeneous Dirichlet boundary conditions, we define a  \ac{fe} lifting $\mathcal{L}_h(g^{\bm{\mu}})$  of the Dirichlet value $g^{\bm{\mu}}$. Let $u_{h,0}^{\bm{\mu}}$ be a \ac{fe} interpolation of $u_0^{\bm{\mu}}$.   
The fully discrete problem at time value $n \in \mathbb{N}(N_t)$  reads: given ${u}_{h,n-1}^{\bm{\mu}} \in \mathcal{V}_h^0$, compute $\forall v_h \in \mathcal{V}_h^0$
\begin{align}
    \label{eq: weak heat equation}
    {u}_{h,n}^{\bm{\mu}} \in \mathcal{V}_{h}^{0} + \mathcal{L}_h(g^{\bm{\mu}})\ : \ \delta^{-1} ({u}_{h,n}^{\bm{\mu}},v_h) +  a(u_{h,n}^{\bm{\mu}},v_h;t_n,\bm{\mu}) = l(v_h;t_n,\bm{\mu}) + \delta^{-1} (u_{h,n-1}^{\bm{\mu}},v_h),
\end{align}
where 
\begin{align}
    \label{eq: forms heat equation}
    a(u,v;t,\bm{\mu}) = \int_\Omega \alpha^{\bm{\mu}}(t)\nabla u \cdot \nabla v \ d\Omega, \qquad 
    l(v;t,\bm{\mu}) = \int_\Omega f^{\bm{\mu}}(t)  v \ d\Omega + \int_{\Gamma_N} h^{\bm{\mu}}(t) v \ d\Gamma. 
\end{align}
We can write this problem in algebraic form as:
\begin{align}
\label{eq: theta method}
    \left(\delta^{-1}\bm{M}_{s,s}+ \bm{A}_{s,s}^{\bm{\mu}}(t_n)\right)\bm{U}^{\bm{\mu}}_{s,n} = 
    \bm{L}_s^{\bm{\mu}}(t_n) +  \delta^{-1}\bm{M}_{s,s} \bm{U}^{\bm{\mu}}_{s,n-1} , \quad n \in \N(N_t), 
\end{align}
where $\bm{U}^{\bm{\mu}}_{s,n} \in \R^{N_s}$ is the vector of \acp{dof} of $u_{h,n}^{\bm{\mu}}$. 
$\bm{M}_{s,s}$ (which is independent of $(t,\bm{\mu})$ throughout this work), $\bm{A}_{s,s}^{\bm{\mu}}$ and $\bm{L}_s^{\bm{\mu}}$ are the mass matrix, the stiffness matrix and the \ac{rhs} vector, respectively. We also introduce the symmetric, positive definite matrix $\bm{X}_{s,s} = \bm{M}_{s,s} + \bm{A}_{s,s}$, representing the $\mathcal{V}$ inner product on $\mathcal{V}_h$. (Note that $\bm{A}_{s,s}$ is the discrete Laplacian with unit diffusivity.) We can also state the space-time algebraic system
\begin{align}
    \label{eq: space time FOM}
    \bm{K}^{\bm{\mu}}_{st,st}\bm{U}^{\bm{\mu}}_{st}=\bm{L}^{\bm{\mu}}_{st},
\end{align}
where $\bm{K}^{\bm{\mu}}_{st,st}$ is a bi-diagonal block-matrix with $N_t$ diagonal blocks $\delta^{-1} \bm{M}_{s,s}+ \bm{A}_{s,s}^{\bm{\mu}}(t_i)$ and lower diagonal blocks $\delta^{-1} \bm{M}_{s,s}$. We can define a global spatio-temporal norm matrix $\bm{X}_{st,st}$, which is a block-diagonal matrix with $N_t$ blocks $\delta \bm{X}_{s,s}$; the $\delta$ serves for this matrix to represent the $L^2(0,T;\mathcal{V})$ product for the \ac{fe} basis defined above. We also introduce the $(\bm{X}_{st,st},\bm{X}_{st,st}^{-1})$ space-time matrix norm: 
\begin{equation}
    \label{eq: X,X^-1 norm}
    \norm{\bm{K}^{\bm{\mu}}_{st,st}}_{\bm{X}_{st,st},\bm{X}_{st,st}^{-1}} =
    \sup\limits_{\bm{V}_{st}} \frac{ \| \bm{K}^{\bm{\mu}}_{st,st} \bm{V}_{st} \|_{\bm{X}_{st,st}^{-1}}}{\norm{\bm{V}_{st}}_{\bm{X}_{st,st}} }. 
\end{equation} 
In the course of this work, we use the property
\begin{equation}
    \label{eq: X,X^-1 norm inequality}
    \norm{\bm{K}^{\bm{\mu}}_{st,st}}_{\bm{X}_{st,st},\bm{X}_{st,st}^{-1}} \leq
    \sup\limits_{\bm{V}_{st}} \frac{\norm{\bm{X}^{-1/2}_{st,st}}_2 \norm{ \bm{K}^{\bm{\mu}}_{st,st} \bm{V}_{st} }_2}{\norm{\bm{V}_{st}}_{\bm{X}_{st,st}}} 
    \leq \norm{\bm{X}^{-1}_{st,st}}_2 \norm{\bm{K}^{\bm{\mu}}_{st,st}}_2, 
\end{equation}
where $\norm{\bm{X}_{st,st}^{-1}}_2$ acts as an equivalence constant between the $(\bm{X}_{st,st},\bm{X}_{st,st}^{-1})$ and $\ell^2$ norms, and scales with $\delta^{-1}(h^{-d}+h^{2-d})$. (In this work, the $\norm{\cdot}_2$ of a square matrix is to be understood as its spectral radius.) We lastly introduce the coercivity constant of $\bm{K}^{\bm{\mu}}_{st,st}$
\begin{align}
    \label{eq: inf sup stability constant}
   \beta =
   \inf\limits_{\bm{V}_{st}}
   \frac{ \| \bm{K}^{\bm{\mu}}_{st,st} \bm{V}_{st} \|_{\bm{X}_{st,st}^{-1}}}{\norm{\bm{V}_{st}}_{\bm{X}_{st,st}} } 
   = \norm{\bm{X}_{st,st}^{-1/2}\bm{K}^{\bm{\mu}}_{st,st}\bm{X}_{st,st}^{-1/2}}_2,
\end{align}
which is positive by virtue of the coercivity of $\bm{K}^{\bm{\mu}}_{st,st}$. When $N_{st} = N_sN_t$ is large, solving \eqref{eq: theta method} for a variety of parameters becomes prohibitively expensive, even when exploiting parallel computations \cite{lassila2013model}. In this scenario, devising a \ac{rom} such as the \ac{strb} method (which we address in the following subsections) that overcomes this issue becomes paramount.

\subsection{Offline phase}
\label{subs: offline phase}
We consider a set of parameters $\mathcal{D}_{\mathrm{off}} = \{\bm{\mu}_k\}_{k=1}^{N_{\bm{\mu}}} \subset \mathcal{D}$, and the \ac{fom} snapshots hypermatrix $\bm{U}_{s,t,\bm{\mu}}$, with $\bm{U}_{s,t}^{\bm{\mu}_k}$ being the solution of \eqref{eq: space time FOM} for $\bm{\mu}_k$. By following the approaches in Sect.~\ref{sec: notation}, we can extract from $\bm{U}_{s,t,\bm{\mu}}$ the spatial basis $\bm{\Phi}_{s, \hat s}$, the temporal basis $\bm{\Phi}_{t, \hat t}$, and consequently the space-time basis $\bm{\Phi}_{st, \hat s\hat t}$.

We now provide the details regarding the choice of the \ac{tpod} ranks $n_s$ and $n_t$. In the following discussions, we take the spatial basis as a point of reference to avoid redundancy. We start by recalling the following standard result of \ac{tpod}. We refer to \cite{quarteroni2015reduced} for additional details. 
\begin{proposition}
    \label{prop:l2_optimality_basis}
    Let $\{\sigma_{s,i}\}_{i=1}^{N_{\sigma_s}}$ be the $N_{\sigma_s} = \min \{N_s,N_tN_{\bm{\mu}}\}$ singular values (sorted by decreasing absolute value) of $\bm{U}_{s,t\bm{\mu}}$. Then:
    \begin{equation}
        \norm{ \bm{U}_{s,t\bm{\mu}} - \bm{\Phi}_{s, \hat s}\bm{\Phi}_{\hat s, s}\bm{U}_{s,t\bm{\mu}} }_F^2  
        = \sum\limits_{i=n_s+1}^{N_{\sigma_s}}\sigma_{s,i}^2.
    \end{equation}  
\end{proposition}

Prop.~\ref{prop:l2_optimality_basis} provides a criterion for the selection of $n_s$: given a user-defined tolerance $\varepsilon \in \R_+$ representing a level of accuracy, we consider the minimum value of $n_s \in \N(N_{\sigma_s})$ such that
\begin{equation}
    \label{eq: truncated POD}
    \sum\limits_{i=1}^{n_s}\sigma_{s,i}^2 \Big / \sum\limits_{i=1}^{N_{\sigma_s}}\sigma_{s,i}^2 \geq 1-\varepsilon^2.
\end{equation}
We essentially require the relative energy (i.e. the squared $\ell^2$ norm) retained by the $n_s$ columns of $\bm{\Phi}_{s, \hat s}$ to be larger than $1-\varepsilon^2$. We represent with $\spod{\cdot}$ the map that, given the snapshot matrix $\bm{U}_{s,t \bm{\mu}}$, returns the \ac{rb} $\bm{\Phi}_{s, \hat s}$ that holds \eqref{eq: truncated POD space} for $n_s$ in \eqref{eq: truncated POD}. {We omit $\varepsilon$, since we will consider the same value in all appearances.} 
Similarly, $\tpod{\cdot}$ represents the time-axis \ac{tpod} that takes $\bm{U}_{t,\hat s \bm{\mu}}$ and computes $\bm{\Phi}_{t, \hat t}$ that holds \eqref{eq: truncated POD time} up to $\varepsilon$, using a temporal version of the criterion \eqref{eq: truncated POD}. For the \ac{rb} method to provide a good approximation, we must ensure that:
\begin{itemize}
    \item the problem we approximate is reducible (see Sect.~\ref{sec: introduction}). If not, the $\spod{\cdot}$ will select $n_s \sim N_s$, i.e. (close to) the entire column space of $\bm{U}_{s,t\bm{\mu}}$, thus preventing us from achieving a meaningful compression; 
    \item $N_{\bm{\mu}}$ is large enough and the choice of the sampling space is appropriate, otherwise the testing error will likely be (considerably) larger than the training error;
    \item $\varepsilon$ is low enough to attain a good enough accuracy, but not so low that $\spod{\cdot}$ selects $n_s \sim N_s$ basis vectors.
\end{itemize}

\begin{remark}
    \label{rmk: cholesky}
    Given a symmetric, positive definite matrix $\bm{Y}_{s,s}$, we can run a modified $\spod{\cdot}$ so that the resulting spatial basis is orthogonal in the $\bm{Y}_{s,s}$-norm instead of the $\ell^2$ norm. By definition, $\bm{Y}_{s,s}$ admits a Cholesky decomposition $\bm{Y}_{s,s} = \bm{H}_{s,s}^T\bm{H}_{s,s}$, with $\bm{H}_{s,s}$ upper-triangular. We compute $\widebar{\bm{\Phi}}_{s, \hat s} = \spod{\bm{H}_{s,s}\bm{U}_{s,t\bm{\mu}}}$, and then $\bm{\Phi}_{s, \hat s}$ is found as $\bm{\Phi}_{s, \hat s} = \bm{H}_{s,s}^{-1}\widebar{\bm{\Phi}}_{s, \hat s}$. We use $\bm{Y}_{s,s} = \bm{X}_{s,s}$ to compute the spatial basis in
    the numerical experiments in Sect.~\ref{sec: section6}. In time, $\ell^2$-orthogonality of $\bm{\Phi}_{t, \hat t}$ is sufficient (i.e. $\bm{Y}_{t,t} = \bm{I}_{t,t}$).
\end{remark}

\begin{remark}
    Reference \cite{quarteroni2015reduced} suggests the implementation of a computationally cheaper $\spod{\cdot}$, compared to the one in \eqref{eq: truncated POD space}. We introduce the smaller $N_{\sigma_s} \times N_{\sigma_s}$ matrix:
    \begin{equation}
    \label{eq: covariance matrix}
        \bm{U}_{\sigma_s,\sigma_s} = 
        \begin{cases}
            \bm{U}_{t\bm{\mu},s}\bm{U}_{s,t\bm{\mu}}, \quad N_{\sigma_s} = N_{\bm{\mu}}N_t \leq N_s, \\
            \bm{U}_{s,t\bm{\mu}}\bm{U}_{t\bm{\mu},s}, \quad N_{\sigma_s} = N_s < N_{\bm{\mu}}N_t.
        \end{cases}
    \end{equation}
    Let us consider the case in which $ N_{\sigma_s} = N_{\bm{\mu}}N_t$ (i.e., $\sigma_s = t\bm{\mu}$ in terms of subindices). Instead of computing the expensive factorization in \eqref{eq: truncated POD space}, we find $\bm{\Phi}_{\sigma_s,\hat s}$ via \ac{tpod} on $\bm{U}_{\sigma_s,\sigma_s}$, and then we recover the correct basis $\bm{\Phi}_{s, \hat s} = \bm{U}_{s, t\bm{\mu}} \bm{\Sigma}_{\sigma_s,\sigma_s}^{-1} \bm{\Phi}_{\sigma_s,\hat s}$, where $\bm{\Sigma}_{\sigma_s,\sigma_s} = \mathrm{diag}(\sigma_{s,1},\hdots,\sigma_{s,N_{\sigma_s}})$. We use this strategy to run every $\spod{\cdot}$ in our numerical experiments; we also employ an analogous approach for $\tpod{\cdot}$.
\end{remark}

At this point, we have discussed the implementation details for constructing the spatial and temporal bases and, by combining them, the space-time basis. Now, we extend the statement of Prop.~\ref{prop:l2_optimality_basis} to derive a similar result of accuracy for \ac{strb}.
\begin{proposition}
\label{prop: rb approximation error} 
  Let $\widehat{\bm{U}}_{st,\bm{\mu}} = \bm{\Phi}_{st, \hat s \hat t}\bm{\Phi}_{\hat s \hat t, st}\bm{U}_{st,\bm{\mu}}$ be the space-time \ac{tpod} approximation in \eqref{eq: truncated POD space-time}, and let $\widetilde{\bm{U}}_{s,t\bm{\mu}} = \bm{\Phi}_{s, \hat s}\bm{\Phi}_{\hat s, s}\bm{U}_{s,t\bm{\mu}}$ be the hypermatrix of space-only approximations. Moreover, let $\{\sigma_{t,i}\}_{i=1}^{N_{\sigma_t}}$ be the $N_{\sigma_t} = \min \{N_t,n_sN_{\bm{\mu}}\}$ singular values (sorted by decreasing absolute value) of $\widetilde{\bm{U}}_{t,\hat s \bm{\mu}}$. The following error bound holds: 
    \begin{equation}
        \label{eq: rb approximation error}
        \norm{\bm{U}_{st,\bm{\mu}} - \widehat{\bm{U}}_{st,\bm{\mu}}}^2_F \leq \sum\limits_{i=n_s+1}^{N_{\sigma_s}}\sigma_{s,i}^2 + \sum\limits_{i=n_t+1}^{N_{\sigma_t}}\sigma_{t,i}^2.
    \end{equation}
    \begin{proof}
        We split the total error as $\bm{U}_{st,\bm{\mu}} - \widehat{\bm{U}}_{st,\bm{\mu}} = (\bm{U}_{st,\bm{\mu}} - \widetilde{\bm{U}}_{st,\bm{\mu}}) + ( \widetilde{\bm{U}}_{st,\bm{\mu}} - \widehat{\bm{U}}_{st,\bm{\mu}})$, where the first term accounts for the space reduction and the second for the time reduction. By virtue of Prop.~\ref{prop:l2_optimality_basis} we have:
        \begin{equation}
            \label{eq: first term}
            \norm{\bm{U}_{st,\bm{\mu}} - \widetilde{\bm{U}}_{st,\bm{\mu}}}^2_F
            = \norm{\bm{U}_{s,t\bm{\mu}} - \widetilde{\bm{U}}_{s,t\bm{\mu}}}^2_F 
            = \sum\limits_{i=n_s+1}^{N_{\sigma_s}}\sigma_{s,i}^2.
        \end{equation}
        Now we focus on the term $\widetilde{\bm{U}}_{s,t}^{\bm{\mu}} - \widehat{\bm{U}}_{s,t}^{\bm{\mu}}$; exploiting Prop.~\ref{prop:l2_optimality_basis}, the definition of $\tpod{\cdot}$, the definition of Frobenius norm, and the property of invariance to cyclic permutations of the trace of a matrix we have:
        \begin{equation}
            \begin{split}
                \label{eq: second term}
                \norm{\widetilde{\bm{U}}_{st,\bm{\mu}} - \widehat{\bm{U}}_{st,\bm{\mu}}}_F^2
                &= \norm{\bm{\Phi}_{s, \hat s}\bm{\Phi}_{\hat s, s}\bm{U}_{s,\bm{\mu},t} - \bm{\Phi}_{s, \hat s}\bm{\Phi}_{\hat s, s}\bm{U}_{s,\bm{\mu},t}\bm{\Phi}_{t, \hat t}\bm{\Phi}_{\hat t, t}}_F^2 
                = \norm{\bm{\Phi}_{s, \hat s}\widetilde{\bm{U}}_{\hat s,\bm{\mu},t}(\bm{I}_{t,t} - \bm{\Phi}_{t, \hat t}\bm{\Phi}_{\hat t, t})}_F^2
                \\& = \mathrm{tr}(\bm{\Phi}_{s, \hat s} \left( \widetilde{\bm{U}}_{\hat s,\bm{\mu},t}(\bm{I}_{t,t} - \bm{\Phi}_{t, \hat t}\bm{\Phi}_{\hat t, t})\right)\left( (\bm{I}_{t,t} - \bm{\Phi}_{t, \hat t}\bm{\Phi}_{\hat t, t})\widetilde{\bm{U}}_{t,\bm{\mu},\hat s}\right)\bm{\Phi}_{\hat s, s})
                \\& = \mathrm{tr}(\bm{\Phi}_{\hat s, s}\bm{\Phi}_{s, \hat s} \left( \widetilde{\bm{U}}_{\hat s,\bm{\mu},t}(\bm{I}_{t,t} - \bm{\Phi}_{t, \hat t}\bm{\Phi}_{\hat t, t})\right)\left( (\bm{I}_{t,t} - \bm{\Phi}_{t, \hat t}\bm{\Phi}_{\hat t, t})\widetilde{\bm{U}}_{t,\bm{\mu},\hat s}\right))
                \\& = \mathrm{tr}(\left( \widetilde{\bm{U}}_{\hat s,\bm{\mu},t}(\bm{I}_{t,t} - \bm{\Phi}_{t, \hat t}\bm{\Phi}_{\hat t, t})\right)\left( (\bm{I}_{t,t} - \bm{\Phi}_{t, \hat t}\bm{\Phi}_{\hat t, t})\widetilde{\bm{U}}_{t,\bm{\mu},\hat s}\right))
                \\& = \norm{\widetilde{\bm{U}}_{\hat s,\bm{\mu},t}(\bm{I}_{t,t} - \bm{\Phi}_{t, \hat t}\bm{\Phi}_{\hat t, t})}_F^2 = \sum\limits_{i=n_t+1}^{N_{\sigma_t}}\sigma_{t,i}^2.
            \end{split}
        \end{equation}
        The statement follows from \eqref{eq: first term}--\eqref{eq: second term}.
    \end{proof}
\end{proposition}
We mention the following corollary of Prop.~\ref{prop: rb approximation error} and \eqref{eq: truncated POD}:
\begin{equation}
    \label{eq: corollary Frobenius norm}
   \norm{\bm{U}_{st,\bm{\mu}} - \widehat{\bm{U}}_{st,\bm{\mu}}}_F^2 \leq \varepsilon \left(\norm{\bm{U}_{s,t\bm{\mu}}}^2_F + \norm{\widetilde{\bm{U}}_{t,\hat s \bm{\mu}}}^2_F\right).
\end{equation}
Although this error bound is not as sharp as the one in \eqref{eq: rb approximation error}, it has the advantage of relating the error to a user-defined quantity.

\subsection{Online phase}
\label{subs: online phase} 
In this phase, two main tasks are carried out. First, the projection of \eqref{eq: space time FOM} onto the space-time reduced subspace. Second, solving the resulting (linear) reduced system of equations. In our work, we consider a Galerkin projection. We refer to \cite{dalsanto2019algebraic} for a detailed review of the more general Petrov-Galerkin projections in an \ac{rb} context. Adopting an algebraic notation, the \ac{strb} problem reads as
\begin{equation}
\label{eq: space time reduced problem}
\text{find} \quad \widehat{\bm{U}}_{\hat s \hat t}^{\bm{\mu}} \quad \text{such that} \quad 
\bm{\Phi}_{\hat s \hat t,st} \left(\bm{L}^{\bm{\mu}}_{st} - \bm{K}^{\bm{\mu}}_{st,st}\bm{\Phi}_{st, \hat s \hat t}\widehat{\bm{U}}_{\hat s \hat t}^{\bm{\mu}}\right)  = \zero{\hat s \hat t}.
\end{equation}
We can also express \eqref{eq: space time reduced problem} as 
$\widehat{\bm{K}}_{\hat s \hat t,\hat s \hat t}^{\bm{\mu}}\widehat{\bm{U}}_{\hat s\hat t}^{\bm{\mu}} = \widehat{\bm{L}}_{\hat s\hat t}^{\bm{\mu}}$, where
\begin{equation}
    \label{eq:  space time reduced LHS-RHS}
    \widehat{\bm{K}}_{\hat s \hat t,\hat s \hat t}^{\bm{\mu}} = \bm{\Phi}_{\hat s \hat t, st} \bm{K}^{\bm{\mu}}_{st,st}\bm{\Phi}_{st, \hat s \hat t}, \qquad
    \widehat{\bm{L}}_{\hat s \hat t}^{\bm{\mu}} = \bm{\Phi}_{\hat s \hat t, st} \bm{L}^{\bm{\mu}}_{st},
\end{equation}
denote the Galerkin compression of the space-time \ac{lhs} and \ac{rhs}, respectively. These quantities are computed following the steps outlined in \cite{Tenderini2022Mueller}.

\section{Empirical interpolation methods for steady and unsteady parameterized PDEs}
\label{sec: section3}
We commence this Sect.~by recalling the properties of \ac{mdeim} applied to a steady problem in Sect.~\ref{subs: MDEIM steady}. Then, in Sects.~\ref{subs: STD-MDEIM}-\ref{subs: FUN-MDEIM}, we propose new \ac{mdeim} architectures for unsteady applications.
 
\subsection{Review of MDEIM for steady PDEs} 
\label{subs: MDEIM steady}
Let us consider a parameterized, steady, nonsingular matrix $\bm{A}_{s,s}^{\bm{\mu}} \in \R^{N_s \times N_s}$. We introduce the vector of nonzero entries $\bm{A}_s^{\bm{\mu}} \in \R^{N_z}$ associated to $\bm{A}_{s,s}^{\bm{\mu}}$ (the two are congruent by isometry). \ac{mdeim} is a data-driven procedure hinging on: 
\begin{itemize}
    \item the construction of a snapshots hypermatrix $\bm{A}_{s,\bm{\mu}} = \left[ \bm{A}_s^{\bm{\mu}_1} ,\ldots, \bm{A}_s^{\bm{\mu}_{N_{\bm{\mu}}}} \right] \in \R^{N_z \times N_{\bm{\mu}}}$;
    \item running a \ac{tpod} on the snapshots hypermatrix to form a basis $\bm{\Phi}_{s, \hat s}^a \in \R^{N_z \times n^a_s}$ spanning a subspace for $\mathrm{range}(\bm{A}_{s,\bm{\mu}})$;
    \item the iterative selection of a list of indices $\{\mathcal{I}_1,\ldots,\mathcal{I}_{n^a_s}\}$, with $\mathcal{I}_i \in \N(N_z) \ \forall i$, which are chosen to minimize an interpolation error over nested subspaces of $\mathrm{range}(\bm{A}_{s,\bm{\mu}})$, in the $\ell^\infty$ norm.
\end{itemize}
We associate with this list the sampling matrix 
\begin{equation}
    \label{eq: sampling matrix}
    \bm{P}_{s, \hat s}^a = \squarebrackets{\bm{e}_s(\mathcal{I}_1) ,\ldots, \bm{e}_s(\mathcal{I}_{n^a_s}) } \in \{0,1\}^{N_z \times n_s^a},
\end{equation}
where $\bm{e}_s(\mathcal{I}) \in \R^{N_z}$ is the $\mathcal{I}^{th}$ vector of the $N_z$-dimensional canonical basis. The implementation of \ac{mdeim} is summarized in Algorithm \eqref{alg: steady MDEIM}. 

\begin{algorithm}[H]
	\begin{algorithmic}[1]
	\State \textbf{Input:} Snapshots hypermatrix $\bm{A}_{s,\bm{\mu}} = \left[ \bm{A}_s^{\bm{\mu}_1} ,\ldots, \bm{A}_s^{\bm{\mu}_{N_{\bm{\mu}}}} \right]$, tolerance $\varepsilon$
	\State \textbf{Output:} Orthonormal basis $\bm{\Phi}_{s, \hat s}^a \in \R^{N_z \times n_s^a}$, sampling matrix $\bm{P}_{s, \hat s}^a\in \{0,1\}^{N_z \times n_s^a}$
	\State Compute $\bm{\Phi}_{s, \hat s}^a = \spod{\bm{A}_{s,\bm{\mu}}}$
	\State Set $\bm{P}_{s, \hat s}^a = \left[ \bm{e}_s(\mathcal{I}_1)  \right]$, where $\mathcal{I}_1 = \argmax \vert \bm{\phi}^a_{s,1} \vert$
	\ForAll{$ k \in \{2,\ldots,{n_s^a}\}$}
    \State Set $\bm{V}_s = [\bm{\Phi}_{s, \hat s}^a]_{:,k}$, $\bm{V}_{s,\hat s} = [\bm{\Phi}_{s, \hat s}^a]_{:,1:k-1}$
	\State Compute residual $\bm{r} = \bm{V}_s - \bm{V}_{s,\hat s} \left( \bm{P}_{\hat s, s}^a \bm{V}_{s,\hat s} \right)^{-1} \bm{P}_{\hat s, s}^a \bm{V}_s$
	\State Update $\bm{P}_{s, \hat s}^a = \left[ \bm{P}_{s, \hat s}^a, \bm{e}_s(\mathcal{I}_k)\right]$, where $\mathcal{I}_k = \argmax \vert \bm{r} \vert$
	\EndFor
	\end{algorithmic}
	\caption{MDEIM algorithm}
	\label{alg: steady MDEIM}
\end{algorithm}	

At the $k^{th}$ iteration, we compute the error between the $k^{th}$ basis vector and its projection on the first $k-1$ columns of the \ac{mdeim} basis, evaluated at the $k-1$ sampling points encoded in $\bm{P}_{s, \hat s}^a$. Then, $\bm{P}_{s, \hat s}^a$ is updated by adding the index corresponding to the largest entry (in absolute value) of the residual. The online phase of \ac{mdeim} consists in computing the approximation
\begin{align}
    \label{eq: MDEIM approx}
    \bm{A}_s^{\bm{\mu}} \approx \widehat{\bm{A}}_s^{\bm{\mu}} = \bm{\Phi}_{s, \hat s}^a\widehat{\bm{A}}_{\hat s}(\bm{\mu}),
\end{align}
where the vector of coefficients $\widehat{\bm{A}}_{\hat s}$ is computed by enforcing the interpolation condition:
\begin{equation}
	\label{eq: MDEIM interpolation}
    \widehat{\bm{A}}_{\hat s}(\bm{\mu}) = \left( \bm{P}_{\hat s, s}^a \bm{\Phi}_{s, \hat s}^a \right)^{-1}\bm{P}_{\hat s, s}^a\bm{A}_s^{\bm{\mu}}. 
\end{equation}
An approximation $\bm{A}_{s,s}^{\bm{\mu}} \approx \widehat{\bm{A}}_{s,s}^{\bm{\mu}}$ on the full set of \ac{hf} \acp{dof} can be recovered by following the steps in \cite{NEGRI2015431}. We remark that the procedure is well posed, since $\bm{P}_{\hat s, s}^a \bm{\Phi}_{s, \hat s}^a$ is nonsingular \cite{chaturantabut2010nonlinear}.

\begin{remark}
    The \ac{mdeim} approximation $\widehat{\bm{A}}_{s,s}^{\bm{\mu}}$ of a symmetric, positive definite matrix $\bm{A}_{s,s}^{\bm{\mu}}$ preserves the symmetry, but not necessarily the positive definiteness. Alternative approaches that bypass this issue include finding the vector of coefficients $\widehat{\bm{A}}_{\hat s}(\bm{\mu})$ subject to a generalized linear constraint that forces $\widehat{\bm{A}}_{s,s}^{\bm{\mu}}$ to be coercive \cite{doi:10.1137/140959602}, or directly performing interpolation on the manifold of symmetric, positive definite matrices \cite{CARLBERG2013623}. Finally, we underline that $\widehat{\bm{A}}_{s,s}^{\bm{\mu}}$ is always non-singular, even in the case of non-symmetric input matrices \cite{NEGRI2015431}. This observation alone supports the derivation of \emph{a posteriori} bounds for the error that the approximation in \eqref{eq: MDEIM approx} commits. 
\end{remark}

\begin{remark}
    \label{rmk: hyper-reduction}
    As shown in \eqref{eq: MDEIM interpolation}, computing $\widehat{\bm{A}}_{\hat s}(\bm{\mu})$ requires assembling the \ac{hf} matrix $\bm{A}_{s,s}^{\bm{\mu}}$ at the indices encoded in $\bm{P}_{s, \hat s}^a$ during the online phase.
    This can be done efficiently by constructing a reduced integration domain, as stated in \cite{NEGRI2015431}. Since, in a \ac{fe} setting, assembling $\bm{A}_{s,s}^{\bm{\mu}}$ requires performing an element-wise numerical integration, we can devise an efficient \ac{fe} assembling procedure that performs the numerical integration only over the elements whose contribution to the entries of $\bm{P}_{\hat s, s}^a \bm{A}_s^{\bm{\mu}}$ is nonzero.
\end{remark}

\subsection{Standard MDEIM for unsteady PDEs}
\label{subs: STD-MDEIM}
In the framework of unsteady parameterized PDEs, the authors of \cite{NEGRI2015431} present an algorithm that progressively builds a global spatial reduced basis for the state variable and one for each of the nonaffinely parameterized \ac{fe} matrices and vectors. However, even though their \emph{simultaneous system approximation and state-space reduction} approach appears to be reasonable regarding memory footprint, its wall time is considerable. This claim becomes even stronger when the reduction occurs in space and time, as in our framework. For this reason, we propose an alternative procedure to separate the system approximation and state reduction (discussed in Sect.~\ref{sec: section2}) steps. For the former, we can devise a \ac{stdmdeim} that conceptually reuses the same steps of a \ac{mdeim} for steady problems. This method is outlined in Algorithm \ref{alg: STD-MDEIM}. 
\begin{algorithm}[H]
	\caption{STD-MDEIM algorithm.}
	\label{alg: STD-MDEIM}
	\begin{algorithmic}[1]
	\State \textbf{Input:} Snapshots hypermatrix $\bm{A}_{s,t,\bm{\mu}}$ such that $[\bm{A}_{s,t,\bm{\mu}}]_{:,j,k} = \bm{A}_s^{\bm{\mu}_k}(t_j)$, tolerance $\varepsilon$
	\State \textbf{Output:} Orthonormal basis $\bm{\Phi}_{s, \hat s}^a$, sampling matrix $\bm{P}_{s, \hat s}^a$
    \State Compute $\bm{\Phi}_{s, \hat s}^a$, $\bm{P}_{s, \hat s}^a$ by running Algorithm \eqref{alg: steady MDEIM} on $\bm{A}_{s,t\bm{\mu}}$, $\varepsilon$
	\end{algorithmic}
\end{algorithm}	
The \ac{stdmdeim} approximation of the nonaffine stiffness matrix reads, for any $\bm{\mu}$:
\begin{align}
    \label{eq: STD-MDEIM approx}
    \bm{A}_s^{\bm{\mu}}(t) \approx \widehat{\bm{A}}_s^{\bm{\mu}}(t) = \bm{\Phi}_{s, \hat s}^a\widehat{\bm{A}}_{\hat s}(t;\bm{\mu}), 
    \quad \text{where} \quad
    \widehat{\bm{A}}_{\hat s}(t;\bm{\mu}) = \left( \bm{P}_{\hat s, s}^a \bm{\Phi}_{s, \hat s}^a \right)^{-1}\bm{P}_{\hat s, s}^a\bm{A}_s^{\bm{\mu}}(t).
\end{align}
This approach aims to find a spatial approximation of the nonaffine matrix at every time step, and as a result, it does not entail a reduction of the time complexity. Thanks to Prop.~\ref{prop:l2_optimality_basis}, we can bound the error of the \ac{stdmdeim} approximation in \eqref{eq: STD-MDEIM approx} as:
\begin{equation}
    \label{eq: STD MDEIM error bound}
    \norm{\bm{A}_{st}^{\bm{\mu}} - \widehat{\bm{A}}_{st}^{\bm{\mu}}}_2 
    \leq \norm{\bm{A}_{s,t\bm{\mu}} - \widehat{\bm{A}}_{s,t\bm{\mu}}}_F
    \leq \vert\vert ( \bm{P}_{\hat s, s}^a \bm{\Phi}_{s, \hat s}^a )^{-1}\vert\vert_F \norm{(\bm{I}_{s,s} - \bm{\Phi}_{s, \hat s}^a\bm{\Phi}_{\hat s, s}^a)\bm{A}_{s,t\bm{\mu}}}_F
    \leq \varepsilon \chi^a_s \norm{\bm{A}_{s,t\bm{\mu}}}_F,
\end{equation}
where $\chi^a_s = \vert\vert ( \bm{P}_{\hat s, s}^a \bm{\Phi}_{s, \hat s}^a )^{-1}\vert\vert_F$. This quantity represents the error committed by the \ac{mdeim} interpolation. Since \ac{stdmdeim} requires computing an $n_s^a$-dimensional vector at $N_t$ time steps, it achieves a dimension reduction proportional to $N_z/n_s^a$. In applications where $N_t$ is large, recovering a reduction in both space and time is desirable. The method outlined in Sect.~\ref{subs: ST-MDEIM} is designed to achieve this property.

\begin{remark}
    Unfortunately, the theoretical lower bound for $\chi^a_s$ is exponential with respect to the reduced dimension $n_s^a$. As a remedy, we can utilize alternative selection strategies of the \ac{mdeim} indices that enjoy error bounds that grow polynomially. We refer to \cite{Drma__2016,otto2019discrete} for more details. In any case, the exponential bound is observed to be quite pessimistic in practice; the \ac{mdeim} compression achieves far better results in most scenarios. We employ the standard selection process of the \ac{mdeim} indices for the numerical experiments in Sect.~\ref{sec: section6}.
\end{remark}

\subsection{Space-Time MDEIM for unsteady PDEs}
\label{subs: ST-MDEIM}
The \ac{stdmdeim} procedure aims to find a suitable subspace for the efficient approximation of $\bm{A}_s^{\bm{\mu}}(t)$ at a given time instant $t$. We propose the \ac{stmdeim} approach to construct a reduced subspace for the space of the stiffness matrices across every time instant and parameter. The procedure is outlined in Algorithm \eqref{alg: ST-MDEIM}. 
\begin{algorithm}[H]
	\caption{ST-MDEIM algorithm}
	\label{alg: ST-MDEIM}
	\begin{algorithmic}[1]
	\State \textbf{Input:} Snapshots hypermatrix $\bm{A}_{s,t,\bm{\mu}}$ such that $[\bm{A}_{s,t,\bm{\mu}}]_{:,j,k} = \bm{A}_s^{\bm{\mu}_k}(t_j)$, tolerance $\varepsilon$
	\State \textbf{Output:} Orthonormal bases $\bm{\Phi}_{s, \hat s}^a$, $\bm{\Phi}_{t, \hat t}^a$, sampling matrices $\bm{P}_{s, \hat s}^a$, $\bm{P}_{t, \hat t}^a$
    \State Compute $\bm{\Phi}_{s, \hat s}^a$, $\bm{P}_{s, \hat s}^a$ by running Algorithm \eqref{alg: steady MDEIM} on $\bm{A}_{s,t\bm{\mu}}$, $\varepsilon$
	\State Compute $\bm{\Phi}_{t, \hat t}^a$, $\bm{P}_{t, \hat t}^a$ by running Algorithm \eqref{alg: steady MDEIM} on $\bm{A}_{t,s\bm{\mu}}$, $\varepsilon$ (with a $\tpod{\cdot}$ instead of a $\spod{\cdot}$)
	\end{algorithmic}
\end{algorithm}	
The \ac{stmdeim} approximation of the nonaffine stiffness matrix is the generalization of \eqref{eq: STD-MDEIM approx} at every time step, and reads for any $\bm{\mu}$ as:
\begin{equation}
    \label{eq: ST-MDEIM approx}
    \bm{A}_{st}^{\bm{\mu}} \approx \widehat{\bm{A}}_{st}^{\bm{\mu}} = 
    \left( \bm{\Phi}_{s, \hat s}^a \otimes \bm{\Phi}_{t, \hat t}^a \right)  \widehat{\bm{A}}_{\hat s \hat t}(\bm{\mu}), 
    \quad \text{where} \quad
    \widehat{\bm{A}}_{\hat s \hat t}(\bm{\mu}) = \left( \bm{P}^a_{\hat s, s}\bm{\Phi}_{s, \hat s}^a \otimes \bm{P}^a_{\hat t, t}\bm{\Phi}_{t, \hat t}^a \right)^{-1} \left( \bm{P}^a_{\hat s, s} \otimes \bm{P}^a_{\hat t, t} \right) \bm{A}_{st}^{\bm{\mu}}.
\end{equation} 
The quantities $\bm{\Phi}_{t, \hat t}^a \in \R^{N_t \times n_t^a}$, $\bm{P}^a_{t, \hat t} = \squarebrackets{\bm{e}_t(\mathcal{t}_1),\ldots,\bm{e}_t(\mathcal{t}_{n_t^a})} \in \{0,1\}^{N_t \times n_t^a}$ are the temporal basis and sampling matrix, and $\widehat{\bm{A}}_{\hat s \hat t}(\bm{\mu}) \in \R^{n_{st}^a}$ is the space-time reduced vector of coefficients, where $n_{st}^a = n_s^an_t^a$. The term $\bm{P}^a_{\hat s, s} \otimes \bm{P}^a_{\hat t, t}$ identifies a space-time reduced integration domain. Thus, by extending the procedure outlined in Remark \ref{rmk: hyper-reduction} to space-time, the \ac{stmdeim} coefficients can be efficiently computed with an online cost scaling with $n_{st}^a$, i.e. with a dimension reduction proportional to $N_zN_t/n_{st}^a$. By virtue of \eqref{eq: corollary Frobenius norm} we can bound the error of the space-time approximation in \eqref{eq: ST-MDEIM approx} as
\begin{equation}
    \label{eq: ST MDEIM error bound}
    \norm{\bm{A}_{st}^{\bm{\mu}} - \widehat{\bm{A}}_{st}^{\bm{\mu}}}_2 
    \leq \varepsilon \chi^a_{st} \sqrt{\norm{\bm{A}_{s,t\bm{\mu}}}^2_F+\norm{\widetilde{\bm{A}}_{t, \hat s \bm{\mu}}}^2_F},
\end{equation}
where $\chi^a_{st} = \norm{\left( \bm{P}^a_{\hat s, s}\bm{\Phi}_{s, \hat s}^a \otimes \bm{P}^a_{\hat t, t}\bm{\Phi}_{t, \hat t}^a \right)^{-1}}_F$, and $\widetilde{\bm{A}}_{\hat s, t \bm{\mu}} = \bm{\Phi}_{\hat s,s}^a\bm{A}_{s,t\bm{\mu}}$. For \ac{stmdeim} to achieve the same accuracy as \ac{stdmdeim}, we must ensure that the manifold of stiffness matrices is reducible in time, as well as in space. Additionally, depending on the nature of the problem at hand, it might be necessary to select a higher value of $N_{\bm{\mu}}$ to achieve a space-time compression at the specified tolerance. We further elaborate on this claim while discussing our numerical results in Sect.~\ref{sec: section6}.

\subsection{Functional MDEIM for unsteady PDEs}
\label{subs: FUN-MDEIM}
Even by employing the hyper-reducing techniques outlined in Remark \ref{rmk: hyper-reduction}, constructing the snapshots hypermatrix $\bm{A}_{s,t,\bm{\mu}} $ used in Algorithms \ref{alg: STD-MDEIM}-\ref{alg: ST-MDEIM} requires a considerable number of operations. To limit the computational burden of this task, we can devise an approach leveraging a functional representation of $\bm{A}_{s,s}^{\bm{\mu}}(t)$. We recall that this quantity is obtained by integrating the bilinear form defined in \eqref{eq: forms heat equation} over the \ac{fe} quadrature points for all combinations of trial and test functions. Instead of computing a \ac{tpod} directly on several snapshots of the matrix, the idea of this method is first to attain a compression of the nonlinearly $(t,\bm{\mu})$-dependent field $\alpha^{\bm{\mu}}(t)$ via \ac{tpod}, similarly to what is done in \cite{https://doi.org/10.1002/fld.2712,doi:10.1137/10081157X,LASSILA20101583,Antil2014}; and then, to run \ac{mdeim} on the matrix snapshots that arise by plugging in the newly obtained $(t,\bm{\mu})$-independent fields in the original form. This step serves to improve the compression achieved by the aforementioned methods, and when applicable limits the additional complexity introduced by geometrical parameterizations, thus overcoming the inefficiency that \cite{NEGRI2015431} points out as a deterrent from functional-based (e.g. EIM or DEIM) hyper-reducing techniques. Additionally, any problem specificity and intrusive changes to the \ac{fe} solver which (D)EIM might suffer --- as mentioned in \cite{NEGRI2015431} --- can be bypassed by an appropriate implementation of the numerical algorithm (see Remark \ref{rmk: gridap}). Lastly, the functional approaches we propose can be easily extended in space-time, thus recovering the same dimensional reduction achieved by \ac{stmdeim}.

Let us consider the hypermatrix $\bm{\alpha}_{s,t,\bm{\mu}} \in \R^{N_q \times N_t \times N_{\bm{\mu}}}$ of evaluations of $\alpha^{\bm{\mu}}(t)$ on the $N_q$ \ac{fe} quadrature points, from which we extract the spatial basis $\bm{\Phi}_{s, \hat s}^{\alpha} = \spod{\bm{\alpha}_{s,t,\bm{\mu}}} \in \R^{N_q \times n_s^\alpha}$. We also introduce a discontinuous Galerkin space $\mathcal{Q}_h$ whose nodes coincide with the $N_q$ quadrature points. Moreover, we introduce the functions $\phi^\alpha_{s,1},\ldots,\phi^\alpha_{s,n_s^\alpha} \in \mathcal{Q}_h$, obtained by interpolating the $n_s^\alpha$ basis vectors of $\bm{\Phi}_{s, \hat s}^{\alpha}$ on $\mathcal{Q}_h$. (We underline that $\phi^\alpha_{s,i} \cong \bm{\phi}^\alpha_{s,i} \ \forall i$, where $\bm{\phi}^\alpha_{s,i}$ is the $i^{th}$ column of $\bm{\Phi}_{s, \hat s}^{\alpha}$, and thus we make no distinction between the two.) Let us now consider the bilinear form
\begin{equation}
    \label{eq: bilinear form abar}
    \widebar{a}(u,v;\phi^\alpha_{s,i}) = \int_\Omega \phi^\alpha_{s,i} \nabla u \cdot \nabla v \ d\Omega \qquad \forall i \in \N(n_s^\alpha), 
\end{equation}
obtained by substituting $\alpha^{\bm{\mu}_k}(t)$ with $\phi^\alpha_{s,i}$ in the bilinear form defined in \eqref{eq: forms heat equation}. We now introduce the \ac{fe} matrix $\widebar{\bm{A}}_{s,s}(\phi^\alpha_{s,i})$ obtained by performing numerical integration on the form in \eqref{eq: bilinear form abar}, the corresponding vector of nonzero entries $\widebar{\bm{A}}_s(\phi^\alpha_{s,i})$ and the matrix that collects all these vectors
\begin{equation}
    \label{eq: snap matrix Abar}
    \widebar{\bm{A}}_{s, \hat s} = \squarebrackets{\widebar{\bm{A}}_s(\phi^\alpha_{s,1}) ,\ldots, \widebar{\bm{A}}_s(\phi^\alpha_{s,n_s^\alpha})}.
\end{equation}
The scope of this \ac{funmdeim} is the computation of $\widebar{\bm{A}}_{s, \hat s}$ itself, from which we can then retrieve a matrix approximation by following the steps outlined in Sect.~\ref{subs: STD-MDEIM}. Specifically, we compute $\widebar{\bm{\Phi}}_s^a$ and $\widebar{\bm{P}}_s^a$ by applying \ac{stdmdeim} on $\widebar{\bm{A}}_{s, \hat s}$, and we approximate $\bm{A}_s^{\bm{\mu}}(t)$ for any $(t,\bm{\mu})$ as:
\begin{align}
    \label{eq: FUN-MDEIM approx}
    \bm{A}_s^{\bm{\mu}}(t) \approx \widehat{\bm{A}}_s^{\bm{\mu}}(t) = \widebar{\bm{\Phi}}_{s, \hat s}^a\widehat{\bm{A}}_{\hat s}(t;\bm{\mu}), 
    \quad \text{where} \quad
    \widehat{\bm{A}}_{\hat s}(t;\bm{\mu}) = \left( \widebar{\bm{P}}_{\hat s, s}^a \widebar{\bm{\Phi}}_{s, \hat s}^a \right)^{-1}\widebar{\bm{P}}_{\hat s, s}^a\bm{A}_s^{\bm{\mu}}(t).
\end{align}
The approximation above is derived without assembling a huge snapshots hypermatrix and running an expensive \ac{tpod}, as is done in \ac{stdmdeim}. Instead, \ac{funmdeim} proposes to compute a small hypermatrix $\bm{\alpha}_{s,t,\bm{\mu}}$ ($N_q \ll N_z$), compress it with a relatively inexpensive \ac{tpod}, and then run \ac{stdmdeim} on the much smaller $\widebar{\bm{A}}_{s, \hat s}$ ($n_s^\alpha \ll N_tN_{\bm{\mu}}$).
Before formally deriving the error bounds for this method, we enunciate the following result. 
\begin{lemma}
    \label{lmm: FUN-MDEIM}
    Let us consider $\widebar{\bm{A}}_{s, \hat s}$ defined in \eqref{eq: snap matrix Abar}. The following result holds:
    \begin{equation}
        \label{eq: lemma statement}
        \forall \bm{V}_s \in \mathrm{range}(\bm{A}_{s,t\bm{\mu}}) \quad \exists \widebar{\bm{V}}_s \in \mathrm{range}(\widebar{\bm{A}}_{s, \hat s})
        \quad \text{such that} \quad
        \norm{\bm{V}_s - \widebar{\bm{V}}_s}_2 \lesssim \varepsilon \norm{\bm{\alpha}_{s,t \bm{\mu}}}_F.
    \end{equation}
    \begin{proof}
        We recall that for any $j \in \N(N_t)$ and $k \in \N(N_{\bm{\mu}})$, we can express $\bm{A}_{s,s}^{\bm{\mu}_k}(t_j)$ as: 
        \begin{equation}
            \label{eq: definition stiffness}
            \bm{A}_{s,s}^{\bm{\mu}_k}(t_j) = T^a(\bm{\alpha}_s^{\bm{\mu}_k}(t_j)), 
        \end{equation}
        where $T^a : \R^{N_q} \to \R^{N_s \times N_s}$ is a continuous operator representing the \ac{fe} assembling procedure. ($T^a$ depends on hyperparameters related to the \ac{fe} assembly, such as the order of quadrature and mesh, not explicitly mentioned in the definition, as the steps in this proof are independent of them.) We can express $\bm{V}_s$ and $\widebar{\bm{V}}_s$ as $\bm{V}_s = \sum\limits_{j,k}V_{jk}\bm{A}_s^{\bm{\mu}_k}(t_j)$ and $\widebar{\bm{V}}_s = \sum\limits_i V_i\widebar{\bm{A}}_s(\phi^\alpha_{s,i})$. Exploiting the continuity of $T^a$, we have:
        \begin{equation}
            \begin{split}
                \norm{\bm{V}_s - \widebar{\bm{V}}_s}_2 &=
                \norm{\sum\limits_{j,k}V_{jk}\bm{A}_s^{\bm{\mu}_k}(t_j) - \sum\limits_i \widebar{V}_i\widebar{\bm{A}}_s(\phi^\alpha_{s,i})}_2    
                = \norm{\sum\limits_{j,k}V_{jk}\bm{A}_{s,s}^{\bm{\mu}_k}(t_j) - \sum\limits_i \widebar{V}_i\widebar{\bm{A}}_{s,s}(\phi^\alpha_{s,i})}_F
                \\&= \norm{\sum\limits_{j,k}T^a(V_{jk}\bm{\alpha}_s^{\bm{\mu}_k}(t_j)) - \sum\limits_i T^a(\widebar{V}_i\bm{\phi}^\alpha_{s,i})}_2
                \lesssim \norm{\sum\limits_{j,k}V_{jk}\bm{\alpha}_s^{\bm{\mu}_k}(t_j) - \sum\limits_i \widebar{V}_i\bm{\phi}^\alpha_{s,i}}_2
            \end{split}
        \end{equation}
        Now, we choose $\widebar{V}_i = \left(\bm{\phi}^\alpha_{s,i}\right)^T\sum\limits_{j,k}V_{jk}\bm{\alpha}_s^{\bm{\mu}_k}(t_j)$, and thus by Prop.~\ref{prop:l2_optimality_basis} and definition of \ac{tpod}
        \begin{equation}
            \label{eq: proof lmm FUN-MDEIM}
            \norm{\bm{V}_s - \widebar{\bm{V}}_s}_2 
            \lesssim \norm{\sum\limits_{j,k}V_{j,k} \Big(\bm{\alpha}_s^{\bm{\mu}_k}(t_j) - \bm{\Phi}^\alpha_{s,\hat s}\bm{\Phi}^\alpha_{\hat s,s}\bm{\alpha}_s^{\bm{\mu}_k}(t_j)\Big)}_2 
            \lesssim \norm{\bm{\alpha}_{s,t\bm{\mu}} - \bm{\Phi}^\alpha_{s,\hat s}\bm{\Phi}^\alpha_{\hat s,s}\bm{\alpha}_{s,t\bm{\mu}}}_F
            \leq \varepsilon \norm{\bm{\alpha}_{s,t \bm{\mu}}}_F,
        \end{equation}
        which concludes the proof.
    \end{proof}
\end{lemma}
Lemma \ref{lmm: FUN-MDEIM} allows us to retrieve the following error bound for the approximation \eqref{eq: FUN-MDEIM approx}:
\begin{align}
    \label{eq: FUN MDEIM error bound}
    \norm{\bm{A}_{st}^{\bm{\mu}} - \widehat{\bm{A}}_{st}^{\bm{\mu}}}_2
    &\leq \ \widebar{\chi}^a_s \norm{\left( \bm{I}_{s, s} - \widebar{\bm{\Phi}}_{s, \hat s}^a\widebar{\bm{\Phi}}_{\hat s, s}^a \right) \bm{A}_{s,t\bm{\mu}} }_F
    \\&\lesssim \widebar{\chi}^a_s \norm{\left( \bm{I}_{s, s} - \widebar{\bm{\Phi}}_{s, \hat s}^a\widebar{\bm{\Phi}}_{\hat s, s}^a \right) \widebar{\bm{A}}_{s,\hat s} }_F 
    + \widebar{\chi}^a_s \varepsilon \norm{ \bm{I}_{s, s} - \widebar{\bm{\Phi}}_{s, \hat s}^a\widebar{\bm{\Phi}}_{\hat s, s}^a }_F \norm{\bm{\alpha}_{s,t \bm{\mu}}}_F  
    \\&\leq \ \widebar{\chi}^a_s \varepsilon (\norm{\widebar{\bm{A}}_{s,\hat s}}_F + \norm{\bm{\alpha}_{s,t\bm{\mu}}}_F),
\end{align}
where $\widebar{\chi}^a_s = \norm{ ( \widebar{\bm{P}}_{\hat s, s}^a \widebar{\bm{\Phi}}_{s, \hat s}^a )^{-1}}_F$. Note that this bound is not necessarily less sharp than the one in \eqref{eq: STD MDEIM error bound}, and indeed our results in Sect.~\ref{sec: section6} suggest that \ac{stdmdeim} and \ac{funmdeim} exhibit the same order in $\varepsilon$.

We can extend the aforementioned concepts to the space-time case. We introduce the time basis $\bm{\Phi}_{t, \hat t}^\alpha \in \R^{N_t \times n_t^\alpha}$ extracted via time-axis \ac{tpod} on $\bm{\alpha}_{t, s \bm{\mu}}$, and its associated sampling matrix $\bm{P}_{t, \hat t}^\alpha$. This method --- \ac{stfunmdeim} --- aims to approximate the stiffness matrix for any $\bm{\mu}$ as:
\begin{equation}
    \label{eq: STFUN-MDEIM approx}
        \bm{A}_{st}^{\bm{\mu}} \approx \widehat{\bm{A}}_{st}^{\bm{\mu}} = 
        \left( \widebar{\bm{\Phi}}_{s, \hat s}^a \otimes \bm{\Phi}_{t, \hat t}^\alpha \right)  \widehat{\bm{A}}_{st}(\bm{\mu}), 
        \quad \text{where} \quad
        \widehat{\bm{A}}_{st}(\bm{\mu}) = \left( \widebar{\bm{P}}_{\hat s, s}^a \widebar{\bm{\Phi}}_{s, \hat s}^a \otimes \bm{P}^\alpha_{\hat t, t} \bm{\Phi}_{t, \hat t}^\alpha \right)^{-1} \roundbrackets{\widebar{\bm{P}}_{\hat s, s}^a \otimes \bm{P}^\alpha_{\hat t, t}} \bm{A}_{st}^{\bm{\mu}}.
\end{equation}
The implementation of \ac{stfunmdeim} is reported in Algorithm \ref{alg: STFUN-MDEIM} (\ac{funmdeim} can be executed by computing only the spatial quantities in this algorithm). To derive the error bounds for \ac{stfunmdeim}, we enunciate the following space-time generalization of Lemma \ref{lmm: FUN-MDEIM}.

\begin{lemma}
    \label{lmm: STFUN-MDEIM}
    Let $\widebar{\bm{A}}_{s, \hat s}$ be defined as in Lemma \ref{lmm: FUN-MDEIM}, let $\widetilde{\bm{\alpha}}_{\hat s, t \bm{\mu}} = \bm{\Phi}_{\hat s, s}^\alpha \bm{\alpha}_{s,t \bm{\mu}}$ and let $\bm{\Phi}_{t, \hat t}^\alpha$ be the time basis extracted from $\widetilde{\bm{\alpha}}_{t, \hat s \bm{\mu}}$. Then
    \begin{equation}
        \label{eq: lemma statement st}
        \begin{split}
            \forall \bm{V}_{st} \in \mathrm{range}(\bm{A}_{st,\bm{\mu}}) \quad &\exists \widebar{\bm{V}}_{st} \in \mathrm{range}(\widebar{\bm{A}}_{s, \hat s} \otimes \bm{\Phi}_{t, \hat t}^\alpha) 
            \quad \text{such that} \\
            \norm{\bm{V}_{st} - \widebar{\bm{V}}_{st}}_2 \lesssim \varepsilon &\sqrt{\norm{\bm{\alpha}_{s,t \bm{\mu}}}^2_F+\norm{\widetilde{\bm{\alpha}}_{t, \hat s \bm{\mu}}}^2_F}.
        \end{split}  
    \end{equation}
    \begin{proof} 
        We define $n_{st}^{\alpha} = n_s^{\alpha}n_t^{\alpha}$, and we associate with any given space-time index $i \in \N(n_{st}^{\alpha})$ a pair of spatial, temporal indices $(i_s,i_t) \in \N(n_s^{\alpha}) \times \N(n_t^{\alpha})$. Let us initially consider the case $\bm{V}_{st} = \bm{A}^{\bm{\mu}_1}_{st}$; proceding as in Lemma \ref{lmm: FUN-MDEIM}, we have that:
        \begin{align}
            \norm{\bm{A}^{\bm{\mu}_1}_{st} - \widebar{\bm{V}}_{st}}_2 
            &= \norm{\bm{A}^{\bm{\mu}_1}_{st} - \sum\limits_{i=1}^{n_{st}^{\alpha}} \widebar{V}_i [ \widebar{\bm{A}}_{s, \hat s} ]_{:,i_s} \otimes [ \bm{\Phi}_{t, \hat t}^\alpha ]_{:,i_t}}_2
            = \norm{\sum\limits_{n=1}^{N_t} \Big(\bm{A}^{\bm{\mu}_1}_s(t_n) - \sum\limits_{i=1}^{n_{st}^{\alpha}} \widebar{V}_i [ \widebar{\bm{A}}_{s, \hat s} ]_{:,i_s}[ \bm{\Phi}_{t, \hat t}^\alpha ]_{n,i_t}\Big)}_2
            \\&\lesssim \norm{\sum\limits_{n=1}^{N_t} \Big(\bm{\alpha}^{\bm{\mu}_1}_s(t_n) - \sum\limits_{i=1}^{n_{st}^{\alpha}} \widebar{V}_i [ \bm{\Phi}_{s, \hat s}^{\alpha} ]_{:,i_s}[ \bm{\Phi}_{t, \hat t}^\alpha ]_{n,i_t}\Big)}_2
            = \norm{\bm{\alpha}^{\bm{\mu}_1}_{st} - \bm{\Phi}_{st, \hat s \hat t}^{\alpha} \widebar{\bm{V}}_{\hat s \hat t} }_2
        \end{align}
        The statement follows if we consider a generic $\bm{V}_{st} = \sum\limits_k V_k\bm{A}^{\bm{\mu}_k}_{st}$, we choose $\widebar{V}_i = \left(\bm{\phi}^\alpha_{st,i}\right)^T\sum\limits_k V_k \bm{\alpha}_{st}^{\bm{\mu}_k}$, we follow a space-time version of \eqref{eq: proof lmm FUN-MDEIM} and we exploit \eqref{eq: corollary Frobenius norm}.
    \end{proof}
\end{lemma}

Lemma \ref{lmm: STFUN-MDEIM} allows us to derive the following error bound for the approximation \eqref{eq: STFUN-MDEIM approx}:
\begin{equation}
    \label{eq: STFUN MDEIM error bound}
    \norm{\bm{A}_{st}^{\bm{\mu}} - \widehat{\bm{A}}_{st}^{\bm{\mu}}}_2 \lesssim
    \widebar{\chi}^a_{st} \varepsilon \Big(\norm{\widebar{\bm{A}}_{s,\hat s} \otimes \bm{\Phi}^\alpha_{t,\hat t}}_F + \sqrt{\norm{\bm{\alpha}_{s,t\bm{\mu}}}^2_F + \norm{\widetilde{\bm{\alpha}}_{t, \hat s \bm{\mu}}}^2_F}\Big),
\end{equation}
where $\widebar{\chi}^a_{st} = \norm{ \left( \widebar{\bm{P}}_{\hat s, s}^a \widebar{\bm{\Phi}}_{s, \hat s}^a \otimes \bm{P}^\alpha_{\hat t, t} \bm{\Phi}_{t, \hat t}^\alpha \right)^{-1}}_F$.  

As we will show in Sect.~\ref{sec: section6}, \ac{funmdeim} and \ac{stfunmdeim} achieve a cheap offline phase compared to their fully algebraic counterparts, thanks to the compression attained by compressing the data encoded in $\bm{\alpha}_{s,t,\bm{\mu}}$. 
\begin{algorithm}[t!]
	\caption{STFUN-MDEIM algorithm.}
	\label{alg: STFUN-MDEIM}
	\begin{algorithmic}[1]
    \State \textbf{Input:} Snapshots hypermatrix $\bm{\alpha}_{s,t,\bm{\mu}}$ such that $[\bm{\alpha}_{s,t,\bm{\mu}}]_{:,j,k} = \bm{\alpha}_s^{\bm{\mu}_k}(t_j)$, tolerance $\varepsilon$
	\State \textbf{Output:} Orthonormal bases $\widebar{\bm{\Phi}}_{s, \hat s}^a$, $\bm{\Phi}_{t, \hat t}^\alpha$, sampling matrices $\widebar{\bm{P}}_{s, \hat s}^a$, $\bm{P}_{t, \hat t}^\alpha$
    \State Compute $\bm{\Phi}_{s, \hat s}^\alpha = \spod{\bm{\alpha}_{s,t \bm{\mu}}}$
    \State Compute $\bm{\Phi}_{t, \hat t}^\alpha = \tpod{\bm{\alpha}_{t, s \bm{\mu}}}$
    \State Compute $\widebar{\bm{A}}_{s, \hat s}$ as defined in \eqref{eq: snap matrix Abar} 
    \State Compute $\widebar{\bm{\Phi}}_{s, \hat s}^a$, $\widebar{\bm{P}}_{s, \hat s}^a$ by running Algorithm \eqref{alg: steady MDEIM} on $\widebar{\bm{A}}_{s, \hat s}$, $\varepsilon$
    \State Compute $\widebar{\bm{P}}_{t, \hat t}^\alpha$ by running the \emph{for} loop in Algorithm \eqref{alg: steady MDEIM} on $\bm{\Phi}_{t, \hat t}^\alpha$, $\varepsilon$
	\end{algorithmic}
\end{algorithm}

\section{Embedding MDEIM in ST-RB: the ST-MDEIM-RB approaches}
\label{sec: section4}
In this section, we discuss the properties of the methods that we obtain by combining the \ac{mdeim} approximation of the nonaffine \ac{hf} quantities with the \ac{strb} procedure. We can identify three main steps in their implementation:
\begin{itemize}
    \item we find the space and time bases $\bm{\Phi}_{s, \hat s}$ and $\bm{\Phi}_{t, \hat t}$ spanning the reduced solution subspace, following the steps outlined in Sect.~\ref{sec: section2};
    \item we perform \ac{mdeim} of the nonaffine operators following one of the four approaches presented in Sect.~\ref{sec: section3}, retrieving an approximate affine decomposition;
    \item thanks to the previous step, we can efficiently compute the space-time Galerkin projection of the problem at hand.
\end{itemize}
In the following subsections, we firstly present the problem formulation of a ST-MDEIM-RB method, and secondly we discuss \emph{a posteriori} theoretical error bounds for this class of \acp{rom}. We make use of the symbol $\dbhat{\cdot}$, which refers to a double compression: we approximate a quantity with \ac{mdeim} (first compression), then we reduce it with \ac{strb} (second compression).

\subsection{Implementation of ST-MDEIM-RB}
\label{subs: ST-MDEIM-RB method}
The \ac{mdeim} strategies outlined in Sect.~\ref{sec: section3} can be leveraged to recover the following approximated affine decompositions: 
\begin{equation}
    \label{eq: MDEIM approximation heat eq}
    \bm{A}_{s,s}^{\bm{\mu}}(t) \approx \widehat{\bm{A}}_{s,s}^{\bm{\mu}}(t) \cong \widehat{\bm{A}}_s^{\bm{\mu}}(t) = \bm{\Phi}_{s, \hat s}^a\widehat{\bm{A}}_{\hat s} \left(t;\bm{\mu}\right); \qquad \bm{L}_s^{\bm{\mu}}(t) \approx \widehat{\bm{L}}_s^{\bm{\mu}}(t) = \bm{\Phi}_{s, \hat s}^l\widehat{\bm{L}}_s\left(t;\bm{\mu}\right); \qquad 
	\forall t \in [0,T].
\end{equation}
We denote with $\widehat{\bm{K}}_{st,st}^{\bm{\mu}}$ the quantity obtained by substituting 
the stiffness matrix $\bm{K}^{\bm{\mu}}_{st,st}$ with its \ac{mdeim} approximation $\widehat{\bm{A}}_{s,s}^{\bm{\mu}}$. The statement of a ST-MDEIM-RB problem applied to our heat equation reads as:
\begin{equation}
\label{eq: ST-MDEIM-RB method}
\widehat{\bm{U}}_{\hat s \hat t}^{\bm{\mu}} \in \mathbb{R}^{n_{st}}
\ : \  \quad 
\bm{\Phi}_{\hat s \hat t, st} \left(\widehat{\bm{L}}_{st}^{\bm{\mu}} - \widehat{\bm{K}}_{st,st}^{\bm{\mu}}\bm{\Phi}_{st, \hat s \hat t}\widehat{\bm{U}}_{\hat s \hat t}^{\bm{\mu}}\right) = \zero{\hat s \hat t},
\end{equation}
or equivalently, $\dbhat{\bm{K}}_{\hat s \hat t, \hat s \hat t}^{\bm{\mu}}\widehat{\bm{U}}_{\hat s \hat t}^{\bm{\mu}} = \dbhat{\bm{L}}_{\hat s \hat t}^{\bm{\mu}}$, 
where the quantities $\dbhat{\bm{K}}_{\hat s \hat t, \hat s \hat t}^{\bm{\mu}}$ and $\dbhat{\bm{L}}_{\hat s \hat t}^{\bm{\mu}}$ are obtained by projecting the \ac{mdeim} quantities $\widehat{\bm{K}}_{st,st}^{\bm{\mu}}$ and $\widehat{\bm{L}}_{st}^{\bm{\mu}}$ onto the space-time \ac{rb} space. We identify four ST-MDEIM-RB methods: 
\begin{itemize}
    \item STD-STRB, when combining STD-MDEIM with \ac{strb};
    \item ST-STRB, when combining ST-MDEIM with \ac{strb};
    \item FUN-STRB, when combining FUN-MDEIM with \ac{strb};
    \item STFUN-STRB, when combining STFUN-MDEIM with \ac{strb}.
\end{itemize}
The error that a ST-MDEIM-RB method commits with respect to the \ac{fom} can be split into two contributions: the first one is due to the \ac{mdeim} approximation of the nonaffine \ac{fe} quantities (system approximation) and the second one arises from the \ac{strb} compression (state approximation). In Sect.~\ref{subs: error estimates}, we perform an \emph{a posteriori} error analysis for ST-MDEIM-RB methods. 

\subsection{Error analysis of ST-MDEIM-RB}
\label{subs: error estimates} 
In this subsection, we derive \emph{a posteriori} error estimates of the proposed ST-MDEIM-RB approaches under the spatio-temporal $\bm{X}_{st,st}$ norm. For the sake of conciseness, we assume that both the Dirichlet and Neumann boundary conditions are homogeneous; however, the following results can be easily generalized to non-homogeneous cases.

\begin{theorem}
    \label{thm: error of STD-RB-MDEIM X norm}
    Let us consider the well-posed problem given by \eqref{eq: space time FOM}, and its well-posed \ac{stmdeimrb} approximation in \eqref{eq: MDEIM approximation heat eq}. The following approximation holds:
    \begin{align}
        \label{eq: error estimate X norm}
        \norm{\bm{U}^{\bm{\mu}}_{st} - \widehat{\bm{U}}_{st}^{\bm{\mu}}}_{\bm{X}_{st,st}} 
        \leq & \ \beta^{-1} \Big(  \chi^l_s \varepsilon \norm{\bm{L}_{s,t\bm{\mu}}}_F \norm{\bm{X}_{st,st}^{-1/2}}_2 +  \chi^a_s \varepsilon \norm{\bm{A}_{s,t\bm{\mu}}}_F \norm{\bm{X}_{st,st}^{-1}}_2\norm{\widehat{\bm{U}}_{st}^{\bm{\mu}}}_{\bm{X}_{st,st}} 
        \\&+ \norm{\widehat{\bm{L}}^{\bm{\mu}}_{st} - \widehat{\bm{K}}^{\bm{\mu}}_{st,st}\widehat{\bm{U}}_{st}^{\bm{\mu}}}_{\bm{X}_{st,st}^{-1}} \Big),
    \end{align}
    where $\widehat{\bm{U}}_{st}^{\bm{\mu}}$ is the approximate solution obtained with \ac{stdmdeim}.
    \begin{proof}
        From the \ac{fom} \eqref{eq: space time FOM}, we have the identities
        \begin{equation}
            \label{eq: initial contribs}
            \begin{split}
                \bm{K}^{\bm{\mu}}_{st,st}(\bm{U}^{\bm{\mu}}_{st} - \widehat{\bm{U}}^{\bm{\mu}}_{st}) 
                &= \bm{L}^{\bm{\mu}}_{st} - \bm{K}^{\bm{\mu}}_{st,st}\widehat{\bm{U}}^{\bm{\mu}}_{st}
                = \bm{L}^{\bm{\mu}}_{st} - \widehat{\bm{L}}^{\bm{\mu}}_{st} + \widehat{\bm{L}}^{\bm{\mu}}_{st} - \widehat{\bm{K}}^{\bm{\mu}}_{st,st}\widehat{\bm{U}}^{\bm{\mu}}_{st} + \widehat{\bm{K}}^{\bm{\mu}}_{st,st}\widehat{\bm{U}}^{\bm{\mu}}_{st} - \bm{K}^{\bm{\mu}}_{st,st}\widehat{\bm{U}}^{\bm{\mu}}_{st}
                \\&= \left(\bm{L}^{\bm{\mu}}_{st} - \widehat{\bm{L}}^{\bm{\mu}}_{st} + \widehat{\bm{K}}^{\bm{\mu}}_{st,st}\widehat{\bm{U}}^{\bm{\mu}}_{st} - \bm{K}^{\bm{\mu}}_{st,st}\widehat{\bm{U}}^{\bm{\mu}}_{st}\right) + \left(\widehat{\bm{L}}^{\bm{\mu}}_{st} - \widehat{\bm{K}}^{\bm{\mu}}_{st,st}\widehat{\bm{U}}^{\bm{\mu}}_{st}\right) = \bm{E}^{\mathrm{M}}_{st} + \bm{E}^{\mathrm{RB}}_{st},
            \end{split}
        \end{equation}
        which we conveniently split into the \ac{mdeim} error contribution $\bm{E}^{\mathrm{M}}_{st}$ and the \ac{strb} error contribution $\bm{E}^{\mathrm{RB}}_{st}$. Left multiplying \eqref{eq: initial contribs} by $\left(\bm{K}^{\bm{\mu}}_{st,st}\right)^{-1}$, taking the $\bm{X}_{st,st}$ norm and using the definition of $\beta$ yields:
        \begin{equation}
            \label{eq: rewrite contribs}
            \begin{split}
                \norm{\bm{U}^{\bm{\mu}}_{st} - \widehat{\bm{U}}_{st}^{\bm{\mu}}}_{\bm{X}_{st,st}} 
                &\leq \norm{\left(\bm{K}^{\bm{\mu}}_{st,st}\right)^{-1} \bm{E}^{\mathrm{M}}_{st}}_{\bm{X}_{st,st}} + \norm{\left(\bm{K}^{\bm{\mu}}_{st,st}\right)^{-1}\bm{E}^{\mathrm{RB}}_{st}}_{\bm{X}_{st,st}} 
                \\&\leq \norm{\bm{X}_{st,st}^{1/2}\left(\bm{K}^{\bm{\mu}}_{st,st}\right)^{-1}\bm{X}_{st,st}^{1/2}}_2 \norm{\bm{X}_{st,st}^{-1/2}\bm{E}^{\mathrm{M}}_{st}}_2 
                \\&\ \ \ \ +  \norm{\bm{X}_{st,st}^{1/2}\left(\bm{K}^{\bm{\mu}}_{st,st}\right)^{-1}\bm{X}_{st,st}^{1/2}}_2 \norm{\bm{X}_{st,st}^{-1/2}\bm{E}^{\mathrm{RB}}_{st}}_2
                \\&= \beta^{-1} \norm{\bm{E}^{\mathrm{M}}_{st}}_{\bm{X}^{-1}_{st,st}} + \beta^{-1} \norm{\bm{E}^{\mathrm{RB}}_{st}}_{\bm{X}^{-1}_{st,st}}.
            \end{split}
        \end{equation}
        The quantity $\norm{\bm{E}^{\mathrm{RB}}_{st}}_{\bm{X}^{-1}_{st,st}}$ corresponds to the last term in the statement \eqref{eq: error estimate X norm}. On the other hand, the \ac{mdeim} contribution can be bounded by recalling \eqref{eq: STD MDEIM error bound}:
        \begin{equation}
            \label{eq: contrib mdeim}
            \begin{split}
                \norm{\bm{E}^{\mathrm{M}}_{st}}_{\bm{X}^{-1}_{st,st}} 
                &\leq \norm{\bm{L}^{\bm{\mu}}_{st} - \widehat{\bm{L}}^{\bm{\mu}}_{st}}_{\bm{X}^{-1}_{st,st}} + \norm{\bm{K}^{\bm{\mu}}_{st,st} - \widehat{\bm{K}}^{\bm{\mu}}_{st,st}}_{\bm{X}_{st,st},\bm{X}_{st,st}^{-1}}\norm{\widehat{\bm{U}}^{\bm{\mu}}_{st}}_{\bm{X}_{st,st}}
                \\&\leq \norm{\bm{X}^{-1/2}_{st,st}}_2 \norm{\bm{L}^{\bm{\mu}}_{st} - \widehat{\bm{L}}^{\bm{\mu}}_{st}}_2 + \norm{\bm{X}^{-1}_{st,st}}_2 \norm{\bm{A}^{\bm{\mu}}_{st} - \widehat{\bm{A}}^{\bm{\mu}}_{st}}_2\norm{\widehat{\bm{U}}^{\bm{\mu}}_{st}}_{\bm{X}_{st,st}}
                \\&\leq \chi^l_s \varepsilon \norm{\bm{L}_{s,t\bm{\mu}}}_F \norm{\bm{X}^{-1/2}_{st,st}}_2 + \chi^a_s \varepsilon \norm{\bm{A}_{s,t\bm{\mu}}}_F \norm{\bm{X}^{-1}_{st,st}}_2\norm{\widehat{\bm{U}}^{\bm{\mu}}_{st}}_{\bm{X}_{st,st}}.
            \end{split}
        \end{equation}
        From \eqref{eq: contrib mdeim}, the quantity $\norm{\bm{E}^{\mathrm{M}}_{st}}_{\bm{X}^{-1}_{st,st}}$ can be bounded by the sum of the first two terms in the statement \eqref{eq: error estimate X norm}, thus concluding the proof.
    \end{proof}
\end{theorem}
When our problem features a reducible system, the error contribution associated with \ac{mdeim} ($\bm{E}^{\mathrm{M}}_{st}$ in the proof) explicitly vanishes as the quality of the approximation improves. On the other hand, the quantity $\bm{E}^{\mathrm{RB}}_{st}$ is a residual-based \emph{a posteriori} \ac{strb} error estimate which also vanishes with $\varepsilon$ when the state variable is reducible. This can be shown by using the following heuristics: when both the \ac{fom} and \ac{strb} formulations are well posed, the \ac{strb} solution will approximate the projection of $\bm{U}_{st}$ on $\mathrm{range}(\bm{\Phi}_{st, \hat s \hat t})$, which we have shown in \eqref{eq: corollary Frobenius norm} commits an approximation error of order $\varepsilon$.

As a corollary of Th.~\ref{thm: error of STD-RB-MDEIM X norm}, we address the case of the remaining \ac{stmdeimrb} methods.
\begin{corollary}
\label{thm: error of STMDEIMRB X norm}
Let us consider the well-posed problem given by \eqref{eq: space time FOM}, and its well-posed approximation in \eqref{eq: MDEIM approximation heat eq}. \\
The following approximation holds for \ac{stmdeim}:
\begin{equation}
    \label{eq: error estimate STMDEIM X norm}
    \begin{split}
        \norm{\bm{U}^{\bm{\mu}}_{st} - \widehat{\bm{U}}_{st}^{\bm{\mu}}}_{\bm{X}_{st,st}} 
        \leq \ & \beta^{-1} \chi^l_{st} \varepsilon \sqrt{\norm{\bm{L}_{s,t\bm{\mu}}}^2_F+\norm{\widetilde{\bm{L}}_{t, \hat s \bm{\mu}}}^2_F} \norm{\bm{X}_{st,st}^{-1/2}}_2 
        \\&+ \beta^{-1} \chi^a_{st} \varepsilon \sqrt{\norm{\bm{A}_{s,t\bm{\mu}}}^2_F+\norm{\widetilde{\bm{A}}_{t, \hat s \bm{\mu}}}^2_F}\norm{\bm{X}_{st,st}^{-1}}_2\norm{\widehat{\bm{U}}_{st}^{\bm{\mu}}}_{\bm{X}_{st,st}} 
        \\&+ \beta^{-1}\norm{\widehat{\bm{L}}^{\bm{\mu}}_{st} - \widehat{\bm{K}}^{\bm{\mu}}_{st,st}\widehat{\bm{U}}_{st}^{\bm{\mu}}}_{\bm{X}_{st,st}^{-1}}.
    \end{split}
\end{equation}
The following approximation holds for \ac{funmdeim}:
\begin{equation}
    \label{eq: error estimate FUNMDEIM X norm}
    \begin{split}
        \norm{\bm{U}^{\bm{\mu}}_{st} - \widehat{\bm{U}}_{st}^{\bm{\mu}}}_{\bm{X}_{st,st}} 
        \lesssim \ & \beta^{-1} \chi^l_s \varepsilon \Big(\norm{\widebar{\bm{L}}_{s,\hat s}}_F + \norm{\bm{f}_{s,t\bm{\mu}}}_F\Big)\norm{\bm{X}_{st,st}^{-1/2}}_2
        \\&+\beta^{-1}\chi^a_s \varepsilon \Big(\norm{\widebar{\bm{A}}_{s,\hat s}}_F + \norm{\bm{\alpha}_{s,t\bm{\mu}}}_F\Big)\norm{\bm{X}_{st,st}^{-1}}_2\norm{\widehat{\bm{U}}_{st}^{\bm{\mu}}}_{\bm{X}_{st,st}} 
        \\&+ \beta^{-1}\norm{\widehat{\bm{L}}^{\bm{\mu}}_{st} - \widehat{\bm{K}}^{\bm{\mu}}_{st,st}\widehat{\bm{U}}_{st}^{\bm{\mu}}}_{\bm{X}_{st,st}^{-1}}.
    \end{split}
\end{equation}
The following approximation holds for \ac{stfunmdeim}:
\begin{equation}
    \label{eq: error estimate STFUNMDEIM X norm}
    \begin{split} 
        \norm{\bm{U}^{\bm{\mu}}_{st} - \widehat{\bm{U}}_{st}^{\bm{\mu}}}_{\bm{X}_{st,st}} 
        \lesssim \ & \beta^{-1} \widebar{\chi}^l_{st} \varepsilon \Big(\norm{\widebar{\bm{L}}_{s,\hat s} \otimes \bm{\Phi}^l_{t,\hat t}}_F + \sqrt{\norm{\bm{f}_{s,t\bm{\mu}}}^2_F + \norm{\widetilde{\bm{f}}_{t, \hat s \bm{\mu}}}^2_F}\Big)\norm{\bm{X}_{st,st}^{-1/2}}_2
        \\&+\beta^{-1}\widebar{\chi}^a_{st} \varepsilon \Big(\norm{\widebar{\bm{A}}_{s,\hat s} \otimes \bm{\Phi}^\alpha_{t,\hat t}}_F + \sqrt{\norm{\bm{\alpha}_{s,t\bm{\mu}}}^2_F + \norm{\widetilde{\bm{\alpha}}_{t, \hat s \bm{\mu}}}^2_F}\Big)\norm{\bm{X}_{st,st}^{-1}}_2\norm{\widehat{\bm{U}}_{st}^{\bm{\mu}}}_{\bm{X}_{st,st}} 
        \\&+ \beta^{-1}\norm{\widehat{\bm{L}}^{\bm{\mu}}_{st} - \widehat{\bm{K}}^{\bm{\mu}}_{st,st}\widehat{\bm{U}}_{st}^{\bm{\mu}}}_{\bm{X}_{st,st}^{-1}}.
    \end{split}
\end{equation}

\begin{proof}
The proof follows the same steps as the one for Thm.~\ref{thm: error of STD-RB-MDEIM X norm}. The only difference comes when we bound the term $\norm{\bm{X}_{st,st}^{-1/2}\bm{E}^{\mathrm{M}}_{st}}_2$, where instead of using \eqref{eq: STD MDEIM error bound} as we do in \eqref{eq: contrib mdeim} we now employ: \eqref{eq: ST MDEIM error bound} for \ac{stmdeim}; \eqref{eq: FUN MDEIM error bound} for \ac{funmdeim}; and \eqref{eq: STFUN MDEIM error bound} for \ac{stfunmdeim}.
\end{proof}
\end{corollary}

\section{Stokes equations}
\label{sec: section5}
In this section, we briefly discuss the application of the ST-MDEIM-RB method to the unsteady, parameterized Stokes equations. We consider an inf-sup stable pair of \ac{fe} spaces for velocity and pressure. By employing the \ac{be} time integrator, we obtain the following system of equations:
\begin{equation}
    \label{eq: theta method nstokes}
    \begin{cases}
        \begin{bmatrix}
            \frac{1}{\delta}\bm{M}_{s,s} + \bm{A}_{s,s}^{\bm{\mu}}(t_n)
            & -\bm{B}_{s,s}^T 
            \\
            \bm{B}_{s,s} &
            \bm{0}_{s,s}
        \end{bmatrix}
        \begin{bmatrix}
            \bm{U}_{s,n}^{\bm{\mu}} \\
            \bm{P}_{s,n}^{\bm{\mu}}
        \end{bmatrix}
        -
        \begin{bmatrix}
            \frac{1}{\delta}\bm{M}_{s,s}\bm{U}^{\bm{\mu}}_{s,n-1} \\
            \bm{0}_{s}
        \end{bmatrix}
        =
        \begin{bmatrix}
            \bm{L}^{\bm{\mu}}(t_n) \\
            \bm{0}_{s}
        \end{bmatrix} 
        \qquad
        \forall n \in \N(N_t). 
    \end{cases} 
\end{equation}
We refer to \cite{salsa2016numerical,volker2016finite} for a detailed exposition of the discretization of the Stokes problem. We represent with  $\bm{U}^{\bm{\mu}}_{s,n}$ and $\bm{P}^{\bm{\mu}}_{s,n}$  the vectors of \ac{dof} values for velocity and pressure at $t_n$, respectively. $\bm{M}_{s,s}$ is the velocity mass matrix, $\bm{A}_{s,s}^{\bm{\mu}}(t_n)$ accounts for the viscous term, $\bm{B}_{s,s}$ is the discrete divergence matrix and its transpose the discrete gradient; $\bm{M}_{s,s}$ and $\bm{B}_{s,s}$ are $(t,\bm{\mu})$-independent in this work.

The \ac{strb} method for this problem features a Galerkin projection onto a cartesian product subspace, i.e. $\mathrm{span}(\bm{\Phi}^u_{st,\hat s \hat t} \times \bm{\Phi}^p_{st,\hat s \hat t})$, where $\bm{\Phi}^u_{st,\hat s \hat t}$ and $\bm{\Phi}^p_{st,\hat s \hat t}$ are the $N_{st}^u \times n_{st}^u$- and $N_{st}^p \times n_{st}^p$-dimensional \ac{strb} bases for velocity and pressure, respectively. To guarantee the \emph{inf-sup} stability of the \ac{strb} framework, it is necessary to augment the reduced basis for the velocity with additional basis vectors --- called supremizers --- both in space and in time. For more details regarding supremizers, we refer to \cite{rozza2005optimization, ballarin2015supremizer, dalsanto2019hyper,pegolotti2021model,Tenderini2022Mueller}. In this work, we construct the supremizers following the ideas in \cite{Tenderini2022Mueller}. The ST-MDEIM-RB problem for the Stokes equation reads as: find $(\widehat{\bm{U}}_{\hat s \hat t}^{\bm{\mu}},\widehat{\bm{P}}_{\hat s \hat t}^{\bm{\mu}})$ such that 
\begin{equation}
    \label{eq: ST-MDEIM-RB stokes}
    \begin{bmatrix}
    \bm{\Phi}_{\hat s \hat t, st}^u\widehat{\bm{K}}_{st,st}^{\bm{\mu}}
    & \bm{\Phi}_{\hat s \hat t, st}^p\bm{B}^T_{st,st}
    \\
    \bm{\Phi}_{\hat s \hat t, st}^u\bm{B}_{st,st} & \zero{\hat s \hat t, st}
    \end{bmatrix}
    \begin{bmatrix}
        \bm{\Phi}_{st, \hat s \hat t}^u\widehat{\bm{U}}_{\hat s \hat t}^{\bm{\mu}} \\ \bm{\Phi}_{st, \hat s \hat t}^p\widehat{\bm{P}}_{\hat s \hat t}^{\bm{\mu}}
    \end{bmatrix} 
    = 
    \begin{bmatrix}
    \bm{\Phi}_{\hat s \hat t, st}^u\widehat{\bm{L}}_{st}^{\bm{\mu}} \\ 
    \zero{\hat s \hat t}
    \end{bmatrix},
\end{equation}
or equivalently, 
\begin{equation}
\begin{bmatrix}
    \dbhat{\bm{K}}_{\hat s \hat t, \hat s \hat t}^{\bm{\mu}}
    & \widehat{\bm{B}^T}_{\hat s \hat t,\hat s \hat t} 
    \\
    \widehat{\bm{B}}_{\hat s \hat t,\hat s \hat t} & \zero{\hat s \hat t,\hat s \hat t}
    \end{bmatrix}
    \begin{bmatrix}
    \widehat{\bm{U}}_{\hat s \hat t}^{\bm{\mu}} \\ \widehat{\bm{P}}_{\hat s \hat t}^{\bm{\mu}}
    \end{bmatrix} 
    = 
    \begin{bmatrix}
    \dbhat{\bm{L}}_{\hat s \hat t}^{\bm{\mu}} \\ 
    \zero{\hat s \hat t}
    \end{bmatrix}.
  \end{equation}
We underline that the blocks $\widehat{\bm{B}^T}_{\hat s \hat t,\hat s \hat t} \in \R^{n^u_{st}\times n^p_{st}}$ and $\widehat{\bm{B}}_{\hat s \hat t,\hat s \hat t} \in \R^{n^p_{st}\times n^u_{st}}$ are obtained without running  a \ac{mdeim} procedure, given the independence of $\bm{B}_{st,st}$ from $(t,\bm{\mu})$. Additionally, we recall that $\widehat{\bm{B}^T}_{\hat s \hat t,\hat s \hat t} \neq \widehat{\bm{B}}_{\hat s \hat t,\hat s \hat t}^T$ (see \cite{Tenderini2022Mueller} for more details).

\section{Numerical results} 
\label{sec: section6} 
In this section, we investigate the numerical properties of the proposed ST-MDEIM-RB methods in the case of the heat and Stokes equations. In particular, we assess the following:
\begin{itemize}
    \item For the online phase, we measure (1) the errors between the \ac{hf} solutions and the ones obtained with ST-MDEIM-RB, and (2) the computational speedup our algorithms achieve compared to the \ac{hf} simulations. As we wish to empirically verify the correctness of the error estimates derived in Sect.~\ref{subs: error estimates}, we make the measurements for several values of the \ac{tpod} tolerance. In particular, we consider $\varepsilon \in \{10^{-i}\}_{i=2}^{4}$.
    \item For the offline phase, we compare the wall time and memory footprint of the \ac{mdeim} approximation of the tangent matrix, for each of the proposed ST-MDEIM-RB strategies. Since the value of $\varepsilon$ has a minimal impact on these quantities, we present the results obtained with $\varepsilon = 10^{-4}$.
\end{itemize}
We run the simulations on a fixed geometry, which is shown in Fig. \ref{fig: geometries}. For both tests, we select a training set according to a uniform distribution on $\mathcal{D}$, i.e. $\mathcal{D}_{\mathrm{off}} = \{\bm{\mu}_k\}_{k=1}^{N_{\bm{\mu}}}$ with $N_{\bm{\mu}} = 80$, and we compute the \ac{fom} solutions for every $\bm{\mu} \in \mathcal{D}_{\mathrm{off}}$. We use $\mathcal{D}_{\mathrm{off}}$ to construct the \ac{strb} approximation, and only the first $30$ parameters of $\mathcal{D}_{\mathrm{off}}$ to obtain the \ac{mdeim} approximation, as empirical evidence suggests that a smaller number of snapshots suffices for this task. We evaluate the online performance of our algorithms on a testing set $\mathcal{D}_{\mathrm{on}} = \{\bm{\mu}_k\}_{k=1}^{N_{\mathrm{on}}} \subset \mathcal{D}$, with $N_{\mathrm{on}} = 10$, for which we also need to compute the corresponding \ac{hf} solutions. To correctly estimate the performances, it is imperative to ensure that $\mathcal{D}_{\mathrm{off}} \cap \mathcal{D}_{\mathrm{on}} = \emptyset$. The accuracy of the methods is evaluated according to the following measures:
\begin{equation}
\label{eq: error definition}
    E^u = \frac{1}{N_{\mathrm{on}}}\overset{N_{\mathrm{on}}}{\underset{i=1}{\sum}}\dfrac{\norm{\widehat{\bm{U}}_{st}^{\bm{\mu}_i} - \bm{U}_{st}^{\bm{\mu}_i}}_{\bm{X}^u_{st,st}}}{\norm{\bm{U}_{st}^{\bm{\mu}_i}}_{\bm{X}^u_{st,st}}};
    \qquad
    E^p = \frac{1}{N_{\mathrm{on}}}\overset{N_{\mathrm{on}}}{\underset{i=1}{\sum}}\dfrac{\norm{\widehat{\bm{P}}_{st}^{\bm{\mu}_i} - \bm{P}_{st}^{\bm{\mu}_i}}_{\bm{X}^p_{st,st}}}{\norm{\bm{P}_{st}^{\bm{\mu}_i}}_{\bm{X}^p_{st,st}}}.
\end{equation}
The norm matrix $\bm{X}^u_{st,st}$ is equal to $\bm{X}_{st,st}$ and $\bm{X}^p_{st,st} \in \R^{N_{st}^p \times N_{st}^p}$ is the matrix representing the discrete $L^2(\Omega)$-$\ell^2$ norm in space-time (which we employ for the pressure in the Stokes test case). Concerning the computational efficiency of the methods, we compute the speedup defined as the ratio between the average wall time of a \ac{fom} simulation and that of a ST-MDEIM-RB simulation. All the numerical tests are run on a local computer with 66$\si{Gb}$ of RAM and an Intel Core i7 running at 3.40$\si{GHz}$. 

\begin{remark}
    \label{rmk: gridap}
    For the generation of the snapshots, we employ Gridap \cite{Badia2020,Verdugo2022}, a library supporting a \ac{fe} approximation of PDEs in the Julia programming language. Among the numerous useful features of Gridap, we mention that this library resorts to a cell-wise representation of fields while integrating the terms of the weak formulation associated with the problem at hand. This allows us to implement with ease our functional approaches, bypassing the implementation issues that \cite{NEGRI2015431} points out as deterrents from functional methods. Indeed, we can simply compress the nonaffinely parameterized fields in a cell-wise manner (i.e. on the quadrature points), and then define an approximated weak formulation where the parameterized fields are replaced with their compressed versions. The numerical solvers that Gridap makes use of can still be used in this scenario. 
\end{remark}

\begin{center}

    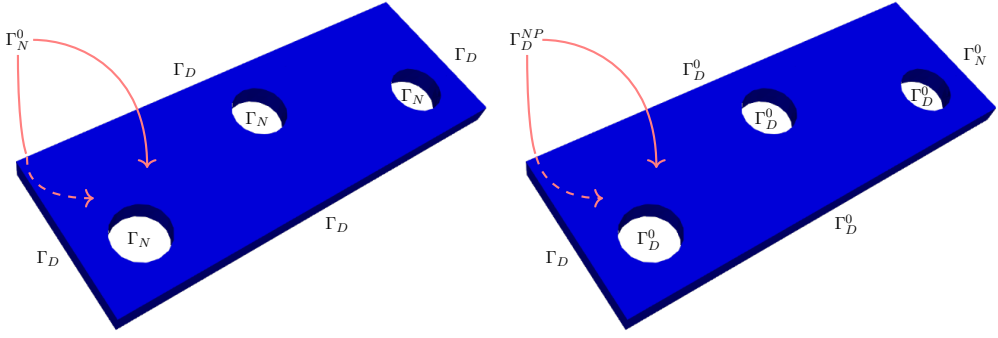
\begin{figure}
        \begin{tikzpicture}
            \node[](image) at (0,0) {\includegraphics[width=\textwidth, height=5cm, keepaspectratio]{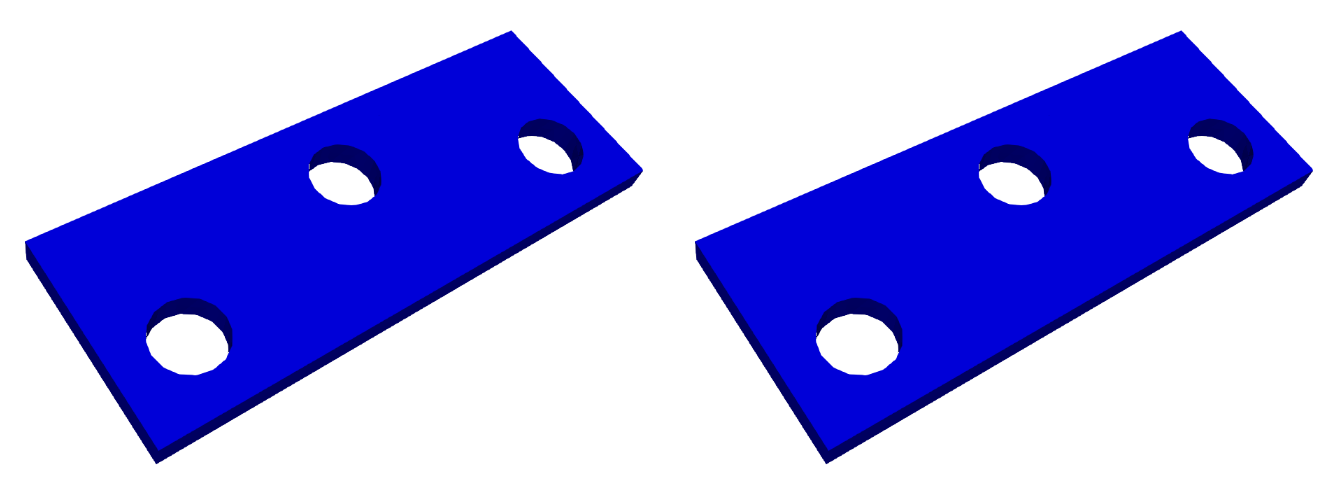}};
            \node[scale=0.6] at (-6,-1.2) {$\Gamma_D$};
            \node[scale=0.6] at (-0.5,1.5) {$\Gamma_D$};
            \node[scale=0.6] at (-2.2,-0.75) {$\Gamma_D$};
            \node[scale=0.6] at (-4.2,1.3) {$\Gamma_D$};
            \node[scale=0.6] at (-4.8,-0.95) {$\Gamma_N$};
            \node[scale=0.6] at (-3.25,0.65) {$\Gamma_N$};
            \node[scale=0.6] at (-1.2,0.95) {$\Gamma_N$};
            \node[scale=0.6] at (-6.4,1.7) {$\Gamma_N^0$};
            \node[scale=0.6] at (0.7,-1.2) {$\Gamma_D$};
            \node[scale=0.6] at (6.2,1.5) {$\Gamma_N^0$};
            \node[scale=0.6] at (2.5,1.3) {$\Gamma_D^0$};
            \node[scale=0.6] at (5.5,0.95) {$\Gamma_D^0$};
            \node[scale=0.6] at (3.45,0.65) {$\Gamma_D^0$};
            \node[scale=0.6] at (1.9,-0.95) {$\Gamma_D^0$};
            \node[scale=0.6] at (4.5,-0.75) {$\Gamma_D^0$};  
            \node[scale=0.6] at (0.3,1.7) {$\Gamma_D^{NP}$};    
            \draw[red!50,->,thick] (-6.2,1.7) .. controls (-5.6,1.7) and (-4.7,1.3) .. (-4.7,0);
            \draw[red!50,thick] (-6.4,1.5) .. controls (-6.4,0.5) and (-6.3,0.2) .. (-6.29,0.19);
            \draw[dashed,red!50,->,thick] (-6.28,0.15) .. controls (-6.27,0.14) and (-6.4,-0.4) .. (-5.4,-0.4);          
            \draw[red!50,->,thick] (0.5,1.7) .. controls (1.1,1.7) and (2,1.3) .. (2,0);
            \draw[red!50,thick] (0.3,1.5) .. controls (0.3,0.5) and (0.4,0.2) .. (0.41,0.19);
            \draw[dashed,red!50,->,thick] (0.42,0.15) .. controls (0.43,0.14) and (0.5,-0.4) .. (1.3,-0.4);
        \end{tikzpicture}
        \vspace{-0.6cm}
        \caption{Geometries used for the heat equation (left) and for the Stokes equations (right). We assign appropriate boundary tags in both cases: $\Gamma_D$, $\Gamma_D^0$, and $\Gamma_D^{NP}$ denote a non-homogeneous, a homogeneous, and a no-penetration Dirichlet boundary; $\Gamma_N$, and $\Gamma_N^0$ denote a non-homogeneous and homogeneous Neumann boundary.}
        \label{fig: geometries}
    \end{figure}
    
    \begin{table}
        \begin{tabular}{ccccccc} 
            \toprule
        
          &\multicolumn{1}{c}{\textbf{Space}}
          &\multicolumn{1}{c}{\textbf{Time}} 
          &\multicolumn{1}{c}{\textbf{Parameters}} 
          &\multicolumn{3}{c}{\textbf{FOM}}\\
            
            \cmidrule(lr){2-2} \cmidrule(lr){3-3} \cmidrule(lr){4-4} \cmidrule(lr){5-7} 
        
          &\multicolumn{1}{c}{$(\bm{L},\bm{H},\bm{W},\bm{R})$}
          &\multicolumn{1}{c}{\bm{$T$}}
          &\multicolumn{1}{c}{\bm{$\mathcal{D}$}}
          &\multicolumn{1}{c}{\textbf{FE Polyn.}}
          &\multicolumn{1}{c}{\bm{$(\bm{N_s^u},\bm{N_s^p},\bm{N_t})$}}
          &\multicolumn{1}{c}{\textbf{WT}} \\
            
            \midrule
        
          \multirow{1}{*}[0.1cm]{}
          &$(4,1.5,0.2,0.25)$ &0.3 &$[1,10]^3$ &Q1 &$(7890,//,60)$ &$17.55 \si{\second}$ \\
        
          \midrule
        
          \multirow{1}{*}[0.15cm]{}
          &$(4,1.5,0.2,0.25)$ &0.15 &$[1,10]^3$ &Q2-P0 &$(93353,5910,60)$ & $798.17 \si{\second}$ \\
            \bottomrule 
        \end{tabular}
        \caption{From left to right: information related to the spatial domain, temporal domain, parameter space, and the \ac{fom}. For the latter, we indicate the \ac{fe} polynomial basis, the spatio-temporal \ac{fom} \acp{dof}, and the average wall time of a \ac{hf} simulation. The quantities $(L,W,H)$ denote the length, width and height of the geometry, whereas $R$ is the radius of the three cylinders.}
        \label{tab: geometries}
    \end{table}
    \end{center}

\subsection{Heat equation}
\label{subs: heat equation results}
In this part, we solve the heat equation introduced in Sect.~\ref{sec: section2}, with the following parametric data:
\begin{equation}
\label{eq: heat equation data}
    \alpha^{\bm{\mu}}(\vec{x},t) = e^{x_1(\sin{t}+\cos{t}) / \sum\limits_i\mu_i};
    \quad
    g^{\bm{\mu}}(\vec{x},t) = \mu_1e^{-\frac{x_1}{\mu_2}}
    \abs{ \sin(\pi t / \mu_3 T) };
    \quad
    h^{\bm{\mu}}(\vec{x},t) = \abs{ \cos(\pi t / \mu_3 T) },
\end{equation}
and $f^{\bm{\mu}}(\vec{x},t) = 1$, $u_0^{\bm{\mu}}(\vec{x}) = 0$. The Dirichlet boundary condition (BC) is imposed strongly, meaning that a lifting operator for the Dirichlet datum is explicitly computed and moved to the \ac{rhs}, as shown in \eqref{eq: weak heat equation}. The unknown temperature is measured in Celsius.

Tb. \ref{tab:heat equation offline} collects the information related to the \ac{mdeim} offline phases. In particular, the functional methods improve by $4$ - $5$ times the performance of their fully algebraic counterparts, both in terms of wall time, and memory allocation. This result reflects the cheaper construction of the snapshot matrices, as well as the reduced cost of the \acp{tpod} we run to extract the bases, as discussed in Sect.~\ref{subs: FUN-MDEIM}. We also infer that the cost of computing the temporal basis of the tangent matrix is marginal, as we can tell from the negligible difference between the \ac{mdeim} methods in space-time, compared to their space-only \ac{mdeim} counterparts. Concerning the online phase, Tb. \ref{tab:heat equation online} presents the efficiency and accuracy measures the ST-MDEIM-RB methods achieve with respect to the \ac{hf} simulations, for three different \ac{tpod} tolerances. We first observe that the error estimates decrease at the same rate as $\varepsilon$, which empirically validates the convergence estimates derived in Sect.~\ref{subs: error estimates}. We deduce that the use of a functional \ac{mdeim} approach does not negatively impact the online accuracy, thus providing empirical evidence of Lemmas \ref{lmm: FUN-MDEIM}-\ref{lmm: STFUN-MDEIM}. We thus infer that the functional \ac{mdeim} retains the same properties as the fully algebraic ones online, while greatly reducing the offline cost, confirming our prior expectations. On the other hand, employing a space-time reducing technique significantly increases the achieved speedup: in particular, ST-STRB and STFUN-STRB are roughly $12$ to $20$ times faster than their space-only counterparts, depending on the value of $\varepsilon$. This is entirely due to the compression of the temporal information these methods achieve. Even greater speedups are expected when considering a larger number of temporal \acp{dof} than in the current setting ($N_t = 60$). Such compression is also responsible for a marginal loss of accuracy compared to the spatial methods, as we can observe from the error estimates. An important takeaway from Tb. \ref{tab:heat equation online} is that it is necessary to select carefully the number of snapshots we employ to construct the \ac{mdeim} bases when employing ST-STRB or STFUN-STRB, if we desire to achieve an accuracy that decays with $\varepsilon$. As we can notice from the error measures in brackets (corresponding to the case in which only $20$ snapshots are used, instead of $30$), the temporal bases fail to correctly approximate the evolution in time of the system when $\varepsilon = 10^{-4}$. We have investigated the source of the error, and we have discovered that there exists some $\bm{\mu}^* \in \mathcal{D}_{on}$ such that $\norm{\widehat{\bm{L}}_{st}^{\bm{\mu}^*} - \bm{L}_{st}^{\bm{\mu}^*}}_2 \eqsim 10\chi^l_s \varepsilon \norm{\bm{L}_{s,t\bm{\mu}}}_F$, i.e. a tenfold increase of \emph{a priori} expectations, thus worsening the average relative error. We conclude that the significant online speedup that characterizes ST-STRB or STFUN-STRB comes at the (potential) cost of performing a more expensive offline phase, because a (potentially) higher number of snapshots is required to ensure an accurate temporal compression. In Fig. \ref{fig:heat equation plots} we compare the temperature obtained by solving the \ac{fom} in \eqref{eq: theta method} with the one computed by employing STFUN-STRB, at $t = 0.15 \si{\second}$, with $\bm{\mu} = \left[7.37,9.58,4.05\right]$, and employing $\varepsilon = 10^{-4}$. (The solutions obtained with the other ST-MDEIM-RB procedures are omitted to avoid redundancy, given that they are graphically indistinguishable from the ones obtained with STFUN-STRB.)

\begin{center}
    \begin{table}[H]
        \begin{tabular}{ccccc} 
        \toprule
        
        &
        \multicolumn{1}{c}{\textbf{STD-STRB}} &
        \multicolumn{1}{c}{\textbf{ST-STRB}} &
        \multicolumn{1}{c}{\textbf{FUN-STRB}} &
        \multicolumn{1}{c}{\textbf{STFUN-STRB}} \\
        
        \midrule
        
        \textbf{WT}
        &\cellcolor{red!25}$51.00 \si{\second}$ 
        &\cellcolor{red!25}$51.04 \si{\second}$ 
        &\cellcolor{blue!25}$11.09 \si{\second}$
        &\cellcolor{blue!25}$11.10 \si{\second}$ \\

        \textbf{MEM}
        &\cellcolor{red!25}$6.26 \si{Gb}$ 
        &\cellcolor{red!25}$6.26 \si{Gb}$ 
        &\cellcolor{blue!25}$1.37 \si{Gb}$
        &\cellcolor{blue!25}$1.37 \si{Gb}$ \\
            
        \bottomrule  
        \end{tabular}
        \caption{Wall time (in $\si{\second}$) and memory footprint (in $\si{Gb}$) of running \ac{mdeim} with $\varepsilon = 10^{-4}$ for the approximation of the tangent matrix of the heat equation.}
        \label{tab:heat equation offline}
    \end{table}
    \begin{table}[H]

        \begin{tabular}{ccccccccc}
        
        \toprule
        
        &
        \multicolumn{2}{c}{\textbf{STD-STRB}} &
        \multicolumn{2}{c}{\textbf{ST-STRB}} &
        \multicolumn{2}{c}{\textbf{FUN-STRB}} &
        \multicolumn{2}{c}{\textbf{STFUN-STRB}} \\
        
        \cmidrule(lr){2-3} \cmidrule(lr){4-5} \cmidrule(lr){6-7}
        \cmidrule(lr){8-9}
        
        $\bm{\varepsilon}$ & 
        \multicolumn{1}{c}{\bm{$SU$}} & 
        \multicolumn{1}{c}{$\bm{E^u / \varepsilon}$} &
        \multicolumn{1}{c}{\bm{$SU$}} & 
        \multicolumn{1}{c}{$\bm{E^u / \varepsilon}$} &
        \multicolumn{1}{c}{\bm{$SU$}} & 
        \multicolumn{1}{c}{$\bm{E^u / \varepsilon}$} &
        \multicolumn{1}{c}{\bm{$SU$}} &
        \multicolumn{1}{c}{$\bm{E^u / \varepsilon}$}\\
        
        \midrule
        
        \multirow{3}{*}{} 
        $\bm{10^{-2}}$ 
        &\cellcolor{red!25} 8.87 &  5.86
        &\cellcolor{blue!25}177.08  & 6.52 $(6.52^*)$
        &\cellcolor{red!25} 10.12 & 5.85
        &\cellcolor{blue!25}180.85 & 6.52 $(6.52^*)$\\
        $\bm{10^{-3}}$ 
        &\cellcolor{red!25} 7.59 & 6.39
        &\cellcolor{blue!25} 113.33 & 6.45 $(6.45^*)$
        &\cellcolor{red!25} 7.33 & 6.40
        &\cellcolor{blue!25} 123.19 & 6.46 $(6.46^*)$\\
        $\bm{10^{-4}}$ 
        &\cellcolor{red!25} 6.30 & 6.78
        &\cellcolor{blue!25} 77.98 & 7.24 $(35.51^*)$
        &\cellcolor{red!25} 6.30 & 6.78
        &\cellcolor{blue!25} 77.98 & 7.24 $(35.51^*)$\\
        
        \bottomrule 
        \end{tabular}
      \caption{Online computational speedup and accuracy achieved by the ST-MDEIM-RB approach compared to the \ac{hf} simulations in the heat equation test case, for three different values of $\varepsilon$. The column relative to the accuracy of ST-STRB and STFUN-STRB displays in brackets the errors obtained when the number of snapshots used to compute the \ac{mdeim} bases is $20$, instead of $30$. The values are not displayed for STD-STRB and FUN-STRB, since changing the number of snapshots from $30$ to $20$ does not affect their error. The results are averaged over 10 different values of $\bm{\mu}$.}
      \label{tab:heat equation online}
    \end{table}
    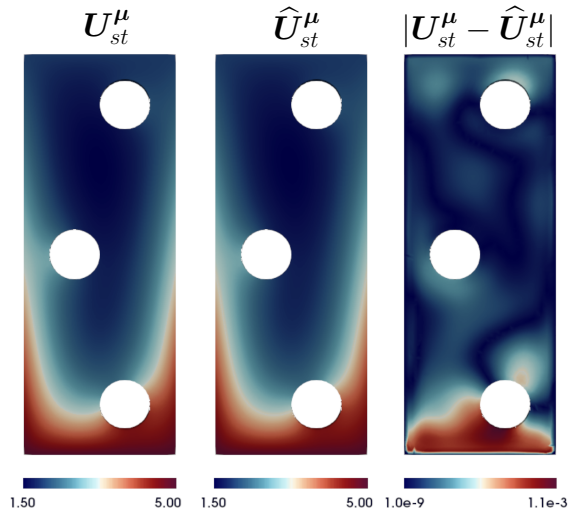
\begin{figure}
        \begin{tikzpicture} 
          \node[](image) at (0,0) {\includegraphics[scale=0.55]{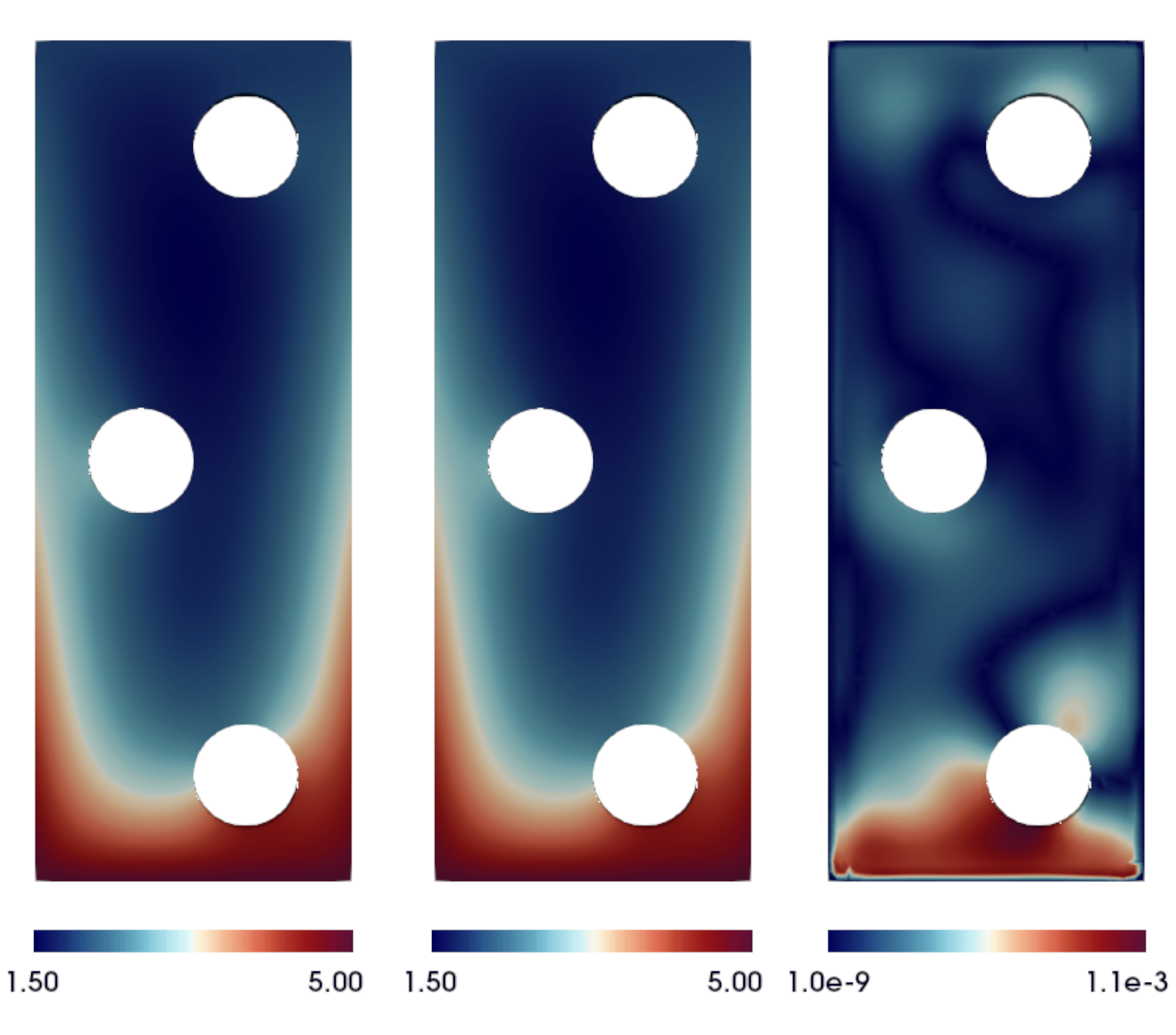}};
          \node[scale=0.6] at (-2.4,3.38) {\huge $\bm{U}_{st}^{\bm{\mu}}$};
          \node[scale=0.6] at (0.1,3.38) {\huge $\widehat{\bm{U}}_{st}^{\bm{\mu}}$};
          \node[scale=0.6] at (2.49,3.38) {\huge $| \bm{U}_{st}^{\bm{\mu}} - \widehat{\bm{U}}_{st}^{\bm{\mu}} |$};
        \end{tikzpicture}
        \vspace{-0.5cm}
        \caption{From left to right: \ac{hf} temperature, STFUN-STRB temperature, and pointwise error. The results are shown from a top view of the geometry, at $t = 0.15 \si{\second}$, with $\bm{\mu} = \left[7.37,9.58,4.05\right]$, and employing $\varepsilon = 10^{-4}$.}
          \label{fig:heat equation plots}	
      \end{figure}
\end{center}

\subsection{Stokes equation}
\label{subs: nstokes equation results}
In this subsection, we solve the Stokes equations with the following parametric data:
\begin{equation*}
\label{eq: nstokes equation data}
    \alpha^{\bm{\mu}}(\vec{x},t) = e^{x_1(\sin{t}+\cos{t}) / \sum\limits_i\mu_i};
    \quad
    \vec{g}^{\bm{\mu}}(\vec{x},t) = -\mu_1 x_2(x_2 - W)\abs{ 1 - \cos(\pi t/T) + \sin(\pi t/\mu_3 T)/\mu_2 }\vec{n}_{\Gamma_D},
\end{equation*}
and $\vec{f}^{\bm{\mu}}(\vec{x},t) = \vec{h}^{\bm{\mu}}(\vec{x},t) = \vec{0}$, $\vec{u}_0^{\bm{\mu}}(\vec{x}) = \vec{0}$. Here, $\vec{n}_{\Gamma_D}$ is the outer normal vector to the inlet $\Gamma_D$, where we impose a parabolic in space, and periodic in time inflow profile. The Dirichlet BC is once again imposed strongly, and we measure the unknowns in the \emph{cgs} (centimeter-gram-second) unit system, i.e. we express the velocity in $cm/s$ and the pressure in $dyn/cm^2$, where $dyn = g \cdot cm / s^2$.

As we show in Tb. \ref{tab: stokes equations offline}, due to the large number of entries of the spatio-temporal tangent matrix, the offline phases of the fully algebraic \ac{mdeim} approaches are too expensive for our computational resources. Therefore, in this Sect.~we are forced only to present the numerical results obtained with the functional methods; these approaches take approximately $1 \si{\minute}$ to compute the \ac{mdeim} approximant of the \ac{lhs} matrix (albeit with a considerable memory footprint). In Tb. \ref{tab: stokes online} we collect the results relative to their online phase. The conclusions we draw are coherent with the ones of the heat equation: both methods achieve an accuracy that decays with $\varepsilon$, and \ac{stfunmdeim} is computationally far more efficient in terms of wall time (depending on the value of $\varepsilon$, it is 23 to 33 times faster than \ac{funmdeim}). The difference in wall time between the two methods is more pronounced than in the heat equation because of the greater number of \acp{dof} characterizing our Stokes equations. In general, we expect this difference in performance to increase with the the number of \acp{dof}, given that the difference in dimension reduction (see Sect.~\ref{subs: STD-MDEIM}-\ref{subs: ST-MDEIM}) will tend to increase with the number of \acp{dof}. The error estimates are smaller than the ones obtained in the heat equation because the state and system approximation errors are both more contained. The former is due to the shorter time interval, which results in a more accurate description of the solution manifold. The latter comes from the simpler parameterization of the problem, as the Neumann BC is homogeneous in this test case. In this test case, extracting the \ac{mdeim} bases from $20$ snapshots instead of $30$ has a marginal effect on the error estimates relative to \ac{stfunmdeim} (the performance worsens when considering $10$-$15$ snapshots). Lastly, in Fig. \ref{fig: stokes} we compare the velocity and pressure obtained by solving the \ac{fom} in \eqref{eq: theta method nstokes} and the ones computed by employing \ac{stfunmdeim}, at $t = 0.075 \si{\second}$, with $\bm{\mu} = \left[5.58,1.20,6.12\right]$, and employing $\varepsilon = 10^{-4}$.

\begin{center}
    \begin{table}[H]
        \begin{tabular}{ccccc} 
            \toprule
            
            &
            \multicolumn{1}{c}{\textbf{STD-STRB}} &
            \multicolumn{1}{c}{\textbf{ST-STRB}} &
            \multicolumn{1}{c}{\textbf{FUN-STRB}} &
            \multicolumn{1}{c}{\textbf{STFUN-STRB}} \\
            
            \midrule
            
            \textbf{WT}
            &\cellcolor{red!25} // 
            &\cellcolor{red!25} //
            &\cellcolor{blue!25}$59.86 \si{\second}$
            &\cellcolor{blue!25}$60.05 \si{\second}$ \\
    
            \textbf{MEM}
            &\cellcolor{red!25} // 
            &\cellcolor{red!25} // 
            &\cellcolor{blue!25}$20.04 \si{Gb}$
            &\cellcolor{blue!25}$20.05 \si{Gb}$ \\
                
            \bottomrule  
            \end{tabular}
        \caption{Wall time (in $\si{\second}$) and memory footprint (in $\si{Gb}$) of running \ac{mdeim} with $\varepsilon = 10^{-4}$ for the approximation of the tangent matrix of the Stokes equations. The results relative to \ac{stdmdeim} and \ac{funmdeim} are not displayed since the memory that they require to perform the task exceeds the maximum memory limit that our computational resources allow.} 
        \label{tab: stokes equations offline}
    \end{table}
    \begin{table}
        \begin{tabular}{cccccc}
        
        \toprule
        
        &
        \multicolumn{2}{c}{\textbf{FUN-STRB}} &
        \multicolumn{2}{c}{\textbf{STFUN-STRB}} & \\
        
        \cmidrule(lr){2-3} \cmidrule(lr){4-5}
        
        $\bm{\varepsilon}$ & 
        \multicolumn{1}{c}{\bm{$SU$}} & 
        \multicolumn{1}{c}{$\{\bm{E^u},\bm{E^p}\}/\bm{ \varepsilon}$} &
        \multicolumn{1}{c}{\bm{$SU$}} &
        \multicolumn{1}{c}{$\{\bm{E^u},\bm{E^p}\}/\bm{ \varepsilon}$}\\
        
        \midrule
        
        \multirow{3}{*}{} 
        $\bm{10^{-2}}$ &\cellcolor{red!25}16.33 &$\{0.36,1.40\}$ &\cellcolor{blue!25}539.31  &$\{1.25,1.32\} (\{3.68,3.36\}^*)$\\
        $\bm{10^{-3}}$ &\cellcolor{red!25}16.26 &$\{0.96,1.39\}$ &\cellcolor{blue!25}411.43 &$\{1.61,2.09\} (\{1.59,2.22\}^*)$\\
        $\bm{10^{-4}}$ &\cellcolor{red!25}16.22 &$\{1.28,1.66\}$ &\cellcolor{blue!25}366.14 &$\{1.31,1.73\} (\{1.34,1.68\}^*)$\\
        
        \bottomrule 
        \end{tabular}
        \caption{Online computational speedup and accuracy achieved by the ST-MDEIM-RB approach with respect to the \ac{hf} simulations in the Stokes equations test case, for three different values of $\varepsilon$. The results are averaged over 10 different values of $\bm{\mu}$.}
        \label{tab: stokes online}
    \end{table}
    \begin{figure}
    \begin{tikzpicture}
        \node[](image) at (0,0) {\includegraphics[scale=0.55]{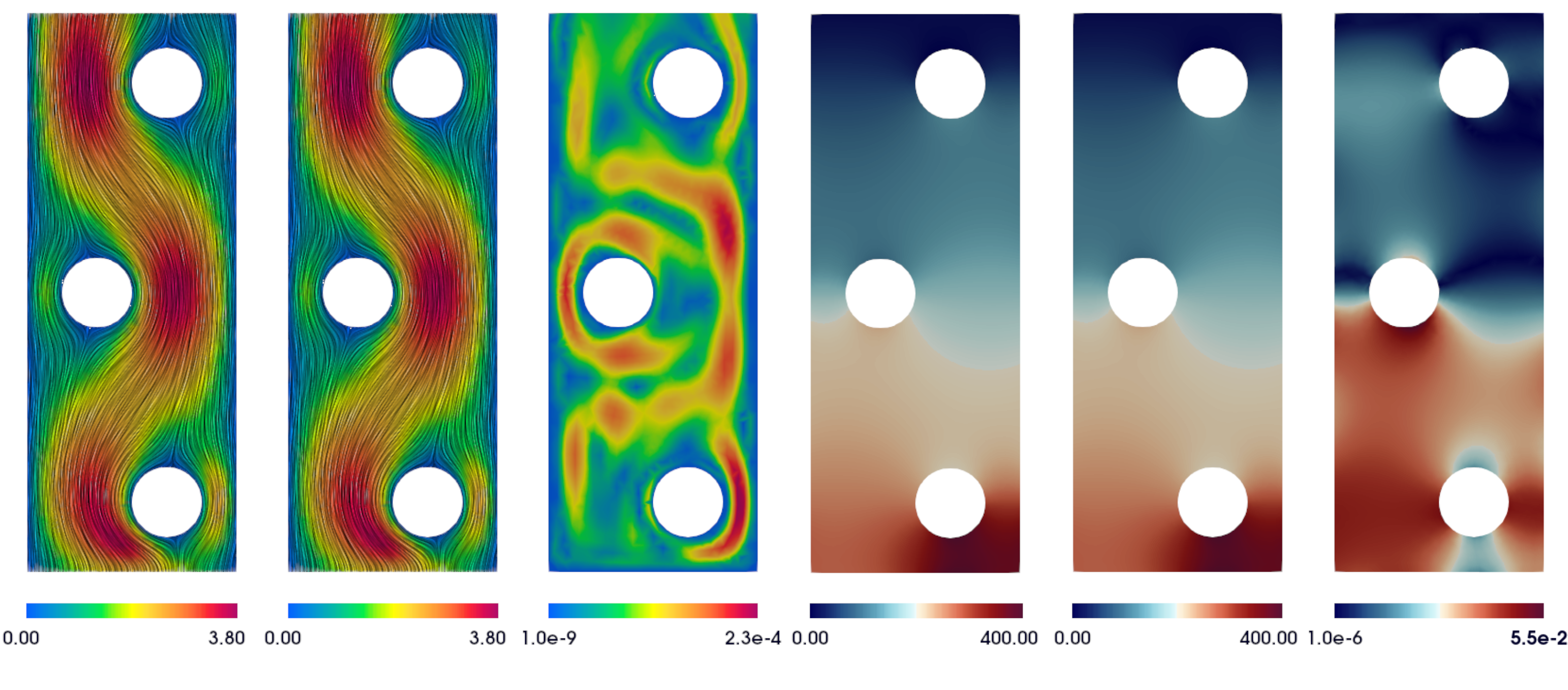}};
        \node[scale=0.6] at (-6.1,3.5) {\huge $\bm{U}_{st}^{\bm{\mu}}$};
        \node[scale=0.6] at (-3.6,3.5) {\huge $\widehat{\bm{U}}_{st}^{\bm{\mu}}$};
        \node[scale=0.6] at (-1.23,3.5) {\huge $| \bm{U}_{st}^{\bm{\mu}} - \widehat{\bm{U}}_{st}^{\bm{\mu}} |$};
        \node[scale=0.6] at (1.23,3.5) {\huge $\bm{P}_{st}^{\bm{\mu}}$};
        \node[scale=0.6] at (3.7,3.5) {\huge $\widehat{\bm{P}}_{st}^{\bm{\mu}}$};
        \node[scale=0.6] at (6.18,3.5) {\huge $| \bm{P}_{st}^{\bm{\mu}} - \widehat{\bm{P}}_{st}^{\bm{\mu}} |$};
    \end{tikzpicture}
    \vspace{-0.5cm}
      \caption{From left to right: surface line integral convolution (SLIC) of the \ac{hf} velocity; SLIC of the \ac{stfunmdeim} velocity; magnitude of the pointwise velocity error; \ac{hf} pressure; \ac{stfunmdeim} pressure; and pointwise pressure error. The results are shown from a top view of the geometry longitudinal mid-section, at $t = 0.075 \si{\second}$, with $\bm{\mu} = \left[5.58,1.20,6.12\right]$ and employing $\varepsilon = 10^{-4}$.} 
      \label{fig: stokes}
    \end{figure}
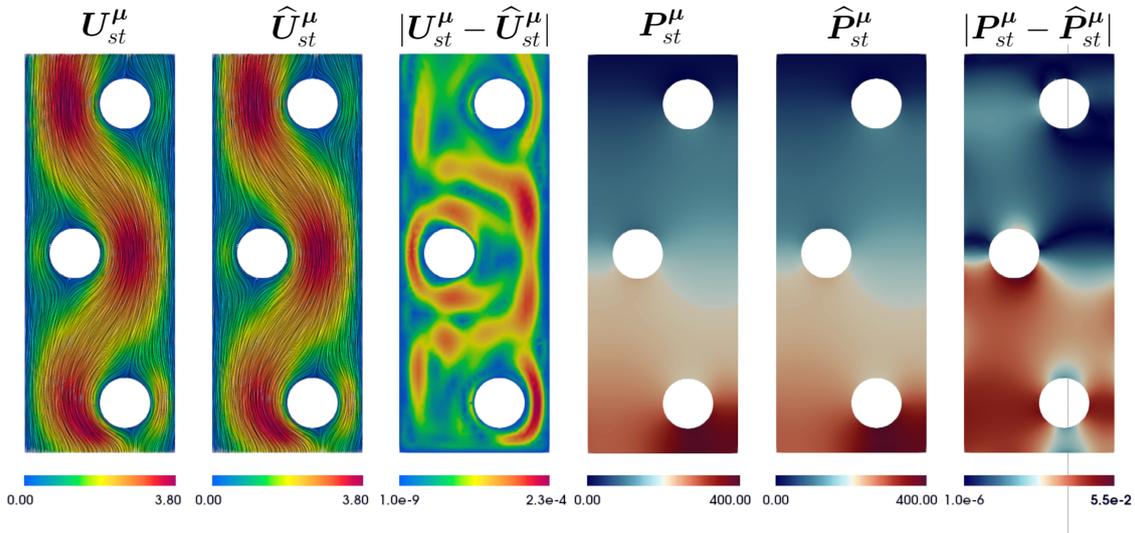
\end{center}

\section{Conclusions} 
\label{sec: section7}
In this work, we discuss the implementation of four ST-MDEIM-RB schemes, consisting of a \ac{mdeim}-based system approximation technique embedded in a \ac{strb} procedure. The first method reduces the spatial dimensionality of nonaffinely parameterized operators at every time step. The second approach achieves an efficient compression of the operator in both space and time. These two strategies proceed in a purely algebraic manner. On the other hand, the remaining two \quotes{functional} methods firstly compress the nonaffine fields evaluated on the quadrature points, thus obtaining a reduced version of the operators; then, they compress the reduced operators in space or space-time following the same steps as their algebraic counterparts. The numerical analysis of the ST-MDEIM-RB algorithms shows that the error decays with a user-defined tolerance that determines the dimension of the \ac{rb} approximation.

We assess the performance of the proposed methods on two test cases, namely a heat equation and a viscous flow governed by the Stokes equations on fixed 3D domains. In each of the tests we consider, the ST-MDEIM-RB algorithms converge as predicted by the theoretical analysis, and substantially shorten the wall time of the \ac{hf} simulations. The methods that rely on a space-time compression are the most efficient ones. On the other hand, the functional approaches are characterized by a cheap offline phase, thanks to the compression of nonaffine fields. Therefore, we can identify a particular approach, i.e. STFUN-STRB, that achieves both cheaper offline and online phases. 

An immediate extension of our work consists of considering nonlinear test cases, such as the Navier-Stokes equations. The proposed methods achieve good results on problems that are reducible when using linear \acp{rom}. We plan to consider space-time nonlinear \acp{rom} in the future, e.g. using neural network-based models with an encoder-decoder structure \cite{choi2019space}.

\printbibliography

\end{document}